\documentclass[a4paper,12pt]{amsart}
\usepackage{amsmath,amsthm,amssymb,latexsym,enumerate,color,hyperref}
\usepackage{graphicx}
\usepackage{geometry}

\usepackage{tikz,fp,ifthen,fullpage} 
\usetikzlibrary{patterns}

\allowdisplaybreaks

\DeclareMathOperator*{\inc}{Incr} 
\begin{document}

\numberwithin{equation}{section}
\numberwithin{figure}{section}

\newtheorem{thm}{Theorem}[section]
\newtheorem{prop}[thm]{Proposition}
\newtheorem{lem}[thm]{Lemma}
\newtheorem{cor}[thm]{Corollary}
\newtheorem{rem}[thm]{Remark}
\newtheorem*{defn}{Definition}

\def\Xint#1{\mathchoice
{\XXint\displaystyle\textstyle{#1}}%
{\XXint\textstyle\scriptstyle{#1}}%
{\XXint\scriptstyle\scriptscriptstyle{#1}}%
{\XXint\scriptscriptstyle\scriptscriptstyle{#1}}%
\!\int}
\def\XXint#1#2#3{{\setbox0=\hbox{$#1{#2#3}{\int}$}
\vcenter{\hbox{$#2#3$}}\kern-.5\wd0}}
\def\ddashint{\Xint=}
\def\dashint{\Xint-}

\newcommand{\DD}{\mathbb{D}}
\newcommand{\NN}{\mathbb{N}}
\newcommand{\ZZ}{\mathbb{Z}}
\newcommand{\QQ}{\mathbb{Q}}
\newcommand{\RR}{\mathbb{R}}
\newcommand{\CC}{\mathbb{C}}
\renewcommand{\SS}{\mathbb{S}}

\newcommand{\ASpace}[1]{\par\vspace*{\stretch{#1}}}
\newcommand\defeq{\mathrel{\stackrel{\rm def}{=}}}
\newcommand{\Frac}[2]{\displaystyle{\frac{#1}{#2}}}
\newcommand{\DS}  {\displaystyle}
\newcommand{\myqed}{\hspace*{\fill}\hbox{\rule[-2pt]{3pt}{6pt}}}
\renewcommand{\theequation}{\arabic{section}.\arabic{equation}}

\newenvironment{prf}[1]{{\noindent \bf Proof of {#1}.}}{\hfill \myqed}

\newcommand{\srom}{\ \rm}
\newcommand{\erom}{\myqed}
\newcommand{\ls}{\{}
\newcommand{\rs}{\}_{n=1,2,\cdots}}
\newcommand{\supp}{\mathop{\mathrm{supp}}}    
\newcommand{\cat}{\mathop{\mathrm{cat}}}                     
\newcommand{\inte}{\mathop{\mathrm{int}}}
\newcommand{\Arcsin}{\mathop{\mathrm{Arcsin}}}
\newcommand{\Arccos}{\mathop{\mathrm{Arccos}}}
\newcommand{\Arctan}{\mathop{\mathrm{Arctan}}}
\newcommand{\BC}{\mathop{\mathrm{BC}}} 
\newcommand{\re}{\mathop{\mathrm{Re}}}   
\newcommand{\im}{\mathop{\mathrm{Im}}}   
\newcommand{\dist}{\mathop{\mathrm{dist}}}  
\newcommand{\link}{\mathop{\circ\kern-.35em -}}
\newcommand{\spn}{\mathop{\mathrm{span}}}   
\newcommand{\ind}{\mathop{\mathrm{ind}}}   
\newcommand{\rank}{\mathop{\mathrm{rank}}}   
\newcommand{\Fix}{\mathop{\mathrm{Fix}}}   
\newcommand{\codim}{\mathop{\mathrm{codim}}}   
\newcommand{\conv}{\mathop{\mathrm{conv}}}   
\newcommand{\epsi}{\mbox{$\varepsilon$}}
\newcommand{\eps}{\mathchoice{\epsi}{\epsi}
{\mbox{\scriptsize\epsi}}{\mbox{\tiny\epsi}}}
\newcommand{\cl}{\overline}
\newcommand{\pa}{\partial}
\newcommand{\ve}{\varepsilon}
\newcommand{\zi}{\zeta}
\newcommand{\Si}{\Sigma}
\newcommand{\cG}{{\mathcal G}}
\newcommand{\cH}{{\mathcal H}}
\newcommand{\cI}{{\mathcal I}}
\newcommand{\cJ}{{\mathcal J}}
\newcommand{\cK}{{\mathcal K}}
\newcommand{\cL}{{\mathcal L}}
\newcommand{\cN}{{\mathcal N}}
\newcommand{\cR}{{\mathcal R}}
\newcommand{\cS}{{\mathcal S}}
\newcommand{\cT}{{\mathcal T}}
\newcommand{\cU}{{\mathcal U}}
\newcommand{\OM}{\Omega}
\newcommand{\B}{\bullet}
\newcommand{\ol}{\overline}
\newcommand{\ul}{\underline}
\newcommand{\vp}{\varphi}
\newcommand{\AC}{\mathop{\mathrm{AC}}}   
\newcommand{\Lip}{\mathop{\mathrm{Lip}}}   
\newcommand{\es}{\mathop{\mathrm{esssup}}}   
\newcommand{\les}{\mathop{\mathrm{les}}}   
\newcommand{\nid}{\noindent}
\newcommand{\pzr}{\phi^0_R}
\newcommand{\pir}{\phi^\infty_R}
\newcommand{\psr}{\phi^*_R}
\newcommand{\pow}{\frac{N}{N-1}}
\newcommand{\ncl}{\mathop{\mathrm{nc-lim}}}   
\newcommand{\nvl}{\mathop{\mathrm{nv-lim}}}  
\newcommand{\la}{\lambda}
\newcommand{\La}{\Lambda}    
\newcommand{\de}{\delta}    
\newcommand{\fhi}{\varphi} 
\newcommand{\ga}{\gamma}    

\newcommand{\core}{\heartsuit}
\newcommand{\diam}{\mathrm{diam}}

\newcommand{\BS}{\PA B_{\ve}(x)}
\newcommand{\lan}{\langle}
\newcommand{\ran}{\rangle}
\newcommand{\tr}{\mathop{\mathrm{tr}}}
\newcommand{\diag}{\mathop{\mathrm{diag}}}
\newcommand{\dv}{\mathop{\mathrm{div}}}
\newcommand{\CI}{\chi_{\{\mu_{p,x}(\ve)\ne u(x+\ve z)\}}(z)}
\newcommand{\tc}{\textcolor}
\newcommand{\TR}{\tc{red}}
\newcommand{\TB}{\tc{blue}}
\newcommand{\TG}{\tc{green}}
\newcommand{\HS}{\phantom{}^t\!z\nabla^2 u(x) z}
\newcommand{\HSE}{\phantom{}^t\!z\nabla^2 u(x+\theta(\ve,z)\ve z) z}
\newcommand{\HSEm}{\phantom{}^t\!z\nabla^2 u(x+\theta\ve_{n_m} z) z}
\newcommand{\HSEn}{\phantom{}^t\!z\nabla^2 u(x+\theta(n,z)\ve z) z}
\newcommand{\HSP}{\phantom{}^t\!z\nabla^2 u(x) z}
\newcommand{\HSEP}{\phantom{}^t\!z\nabla^2 u(x+\theta(\ve,z)\ve z) z}
\newcommand{\HSEPm}{\phantom{}^t\!z\nabla^2 u(x+\ve z) z}
\newcommand{\al}{\alpha}
\newcommand{\be}{\beta}
\newcommand{\Om}{\Omega}
\newcommand{\na}{\nabla}
\newcommand{\svm}{\mathrm{SV}_{mpr}}
\newcommand{\svh}{\mathrm{SV}^1_{hr}}
\newcommand{\svr}{\mathrm{SV}^2_{hr}}
\newcommand{\cC}{\mathcal{C}}
\newcommand{\cM}{\mathcal{M}}
\newcommand{\nr}{\Vert}
\newcommand{\De}{\Delta}
\newcommand{\cX}{\mathcal{X}}
\newcommand{\cP}{\mathcal{P}}
\newcommand{\om}{\omega}
\newcommand{\si}{\sigma}
\newcommand{\te}{\theta}
\newcommand{\Ga}{\Gamma}

\title{An introduction to the study of critical points \\ of solutions of elliptic and parabolic equations}

\author[R. Magnanini]{R. Magnanini}
    \address{Dipartimento di Matematica ed Informatica ``U.~Dini'',
Universit\` a di Firenze, viale Morgagni 67/A, 50134 Firenze, Italy.}
    \email{magnanin@math.unifi.it}
    \urladdr{http://web.math.unifi.it/users/magnanin}

\begin{abstract}
We give a survey at an introductory level of old and recent results in the study of critical points of solutions of
elliptic and parabolic partial differential equations. To keep the presentation simple, we mainly consider four exemplary boundary value problems: the Dirichlet problem for the Laplace's equation; the torsional creep problem; the case of Dirichlet eigenfunctions for the Laplace's equation; the initial-boundary value problem for the heat equation. We shall mostly address three issues: the estimation of the
local size of the critical set; the dependence of the number of critical points on the boundary values and the geometry of the domain; the location of critical points in the domain.
\end{abstract}


\maketitle

\setcounter{equation}{0}
\setcounter{figure}{0}

\section{Introduction} 

Let $\Om$ be a domain in the Euclidean space $\RR^N$, $\Ga$ be its boundary and $u:\Om\to \RR$  be a differentiable function.  A {\it critical point} of $u$ is a point in $\Om$ at which the gradient $\na u$ of $u$ is the zero vector.  The importance of critical points is evident. At an elementary level, they help us to visualize the graph of $u$, since they are some of its notable points (they are local maximum, minimum, or inflection/saddle points 
of $u$). At a more sophisticated level, if we interpret $u$ and $\na u$ as a gravitational, electrostatic or velocity potential and its underlying  field of force or flow, the critical points are the positions of equilibrium for the field of force or stagnation points for the flow and give information on the topology of the equipotential lines or of the curves of steepest descent (or stream lines) related to $u$.
\par
A merely differentiable function can be very complicated. For instance, Whitney \cite{Wh} constructed a non-constant function of class $C^1$ on the plane with a connected set of critical values (the images of critical points). If we allow enough smoothness, this is no longer possible as Morse-Sard's lemma informs us: indeed, if $u$ is at least of class $C^N$, the set of its critical values must have zero Lebesgue measure and hence the regular values of $u$ must be dense in the image of $u$ (see \cite{AR} for a proof).
\par
When the function $u$ is the solution of some partial differential equation, the situation improves. 
In this survey, we shall consider the four archetypical equations:
\begin{equation*}
\De u=0, \quad  \De u=-1, \quad \De u+\la u=0, \quad u_t=\De u,
\end{equation*}
that is the {\it Laplace's} equation, the {\it torsional creep} equation, the {\it eigenfunction} equation and the {\it heat} equation. 
\par
It should be noticed at this point some important differences between the first and the remaining  
three equations. 
\par
One is that the critical points of harmonic functions --- the solutions of the Laplace's equation --- 
are always ``saddle points'' as it is suggested by the maximum and minimum principles and the fact that
$\De u$ is the sum of the eigenvalues of the hessian matrix $\na^2 u$. 
The other three equations instead admit solutions with maximum or minimum points.  
\par
Also, we know that the critical points of a non-constant harmonic function $u$ on an open set of $\RR^2$ are isolated and can be assigned a sort of finite multiplicity, for they are the zeroes of the holomorphic
function $f=u_x-i u_y$. By means of the theory of quasi-conformal mappings and generalized analytic functions, this result can be extended to solutions of the elliptic equation
\begin{equation}
\label{elliptic}
(a\,u_x+b\,u_y)_x+(b\,u_x+c\,u_y)_y+d\,u_x+e\,u_y=0 
\end{equation}
(with suitable smoothness assumptions on the coefficients) or even to weak solutions of the an elliptic equation in divergence form,
\begin{equation}
\label{elliptic2}
(a\,u_x+b\,u_y)_x+(b\,u_x+c\,u_y)_y=0 \ \mbox{ in } \ \Om,
\end{equation}
even allowing discontinuous coefficients. 
\par
Instead, solutions of the other three equations can show curves
of critical points in $\RR^2$, as one can be persuaded by looking at the solution of the torsional creep equation in a circular annulus with zero boundary values.  
\par
These discrepancies extend to any dimension $N\ge 2$, in the sense that 
it has been shown that the set of the critical points of a non-constant harmonic function (or of a solution of an elliptic equation with smooth coefficients modeled on the Laplace equation) has 
at most locally finite $(N-2)$-dimensional Hausdorff measure, while solutions of equations fashioned on the other three equations have at most locally finite $(N-1)$-dimensional Hausdorff measure.
\par
Further assumptions on solutions of a partial differential equation, such as their behaviour on the boundary and the shape of the boundary itself, can give more detailed information on the number and location of critical points. In these notes, we shall consider the case of harmonic functions with various boundary behaviors and the solutions $\tau$, $\phi$ and $h$ of the following three problems:
\begin{equation}
\label{torsion}
-\De \tau=1 \ \mbox{ in } \ \Om, \quad \tau=0 \ \mbox{ on } \  \Ga;
\end{equation}
\begin{equation}
\label{eigenfunction}
\De \phi+\la \phi=0 
\ \mbox{ in } \ \Om, \quad 
\phi=0 \ \mbox{ on } \ \Ga;
\end{equation}
\begin{eqnarray}
&h_t=\De h \ \mbox{ in } \ \Om\times(0,\infty), \label{heat1}\\
&h=0 \ \mbox{ on } \ \Ga\times(0,\infty), \quad
h=\fhi \ \mbox{ on } \ \Om\times\{ 0\}, \label{heat2}
\end{eqnarray}
where $\fhi$ is a given function. We will refer to \eqref{torsion}, \eqref{eigenfunction}, \eqref{heat1}-\eqref{heat2}, as the {\it torsional creep problem}, the {\it Dirichlet eigenvalue problem},
and the {\it initial-boundary value problem for the heat equation}, respectively.
\par
A typical situation is that considered in  Theorem \ref{th:am1b}: a harmonic function $u$ on a planar domain $\Om$ is given together with a vector field $\ell$ on $\Ga$ of assigned topological degree $D$; the number of critical points in $\Om$ then is bounded in terms of $D$, the Euler characteristic of $\Om$ and the number of proper connected components of the set 
$\{ z\in\Ga: \ell(z)\cdot\na u(z)>0\}$ (see Theorem \ref{th:am1b} for the exact statement).
We shall also see how  this type of theorem has recently been extended to obtain a bound for the number 
of critical points of the Li-Tam Green's function of a non-compact Riemanniann surface of finite type in terms of its genus and the number of its ends.
\par
Owing to the theory of quasi-conformal mappings, Theorem \ref{th:am1b} can be extended to solutions of quite general elliptic equations and, thanks to the work of G. Alessandrini and co-authors, has found effective applications to the study of inverse problems that have as a common denominator the reconstruction of the coefficients of an elliptic equation in a domain from measurements on the boundary of a set of its solutions.
\par
A paradigmatic example is that of Electric Impedence Tomography (EIT) in which a conductivity $\ga$
is reconstructed, as the coefficient of the elliptic equation
$$
\dv(\ga\na u)=0 \ \mbox{ in } \ \Om,
$$
from the so-called Neumann-to-Dirichlet (or Dirichlet-to-Neumann) operator on $\Ga$.
In physical terms, an electrical current (represented by the co-normal derivative
$\ga u_\nu$) is applied on $\Ga$  generating a potential $u$, that is measured on $\Ga$ within a certain error. One wants to reconstruct the conductivity $\ga$ from some of these measurements.
Roughly speaking, one has to solve for the unknown $\ga$ the first order differential equation
$$
\na u\cdot\na\ga+(\De u)\,\ga=0 \ \mbox{ in } \ \Om,
$$ 
once the information about $u$ has been extended from $\Ga$ to $\Om$. 
It is clear that such an equation is singular at the critical points of $u$. Thus, it is helpful to know
{\it a priori} that $\na u$ does not vanish and this can be done via (appropriate generalizations of) Theorem \ref{th:am1b} by choosing
suitable currents on $\Ga$.
\par
The possible presence of maximum and/or minimum points for the solutions of \eqref{torsion}, \eqref{eigenfunction}, or \eqref{heat1}-\eqref{heat2} makes the search for an estimate of the number of critical points a difficult task (even in the planar case). In fact, the mere topological information only results in an estimate of the
{\it signed sum} of the critical points, the sign depending on whether the relevant critical point is an extremal or saddle point. For example, for the solution of \eqref{torsion} or \eqref{eigenfunction}, we only know that the difference between the number of its (isolated) maximum and saddle points (minimum points are not allowed) must equal $\chi(\Om)$, the {\it Euler characteristic} of $\Om$ --- a Morse-type theorem.
Thus, further assumptions, such as geometric information on $\Om$, are needed. More information is also necessary even if we consider the case of harmonic functions in dimension $N\ge 3$.
\par
In the author's knowledge, results on the number of critical points of solutions of \eqref{torsion}, \eqref{eigenfunction}, or \eqref{heat1}-\eqref{heat2} reduce to deduction that their solutions admit
a {\it unique} critical point if $\Om$ is {\it convex}. Moreover, the proof of such results is somewhat indirect: 
the solution is shown to be {\it quasi-concave} --- indeed, log-concave for the cases of  \eqref{eigenfunction} and \eqref{heat1}-\eqref{heat2}, and $1/2$-concave for the case \eqref{torsion} --- and then its analyticity completes the argument.  Estimates of the number of critical points when the domain $\Om$ has more complex geometries would be a significant advance. In this survey, we will
propose and justify some conjectures.
\par
The problem of locating critical points is also an interesting issue. The first work on this subject dates back to Gauss \cite{Ga}, who proved that the critical points of a complex polynomial are its, if they are multiple, and the equilibrium points of the gravitational field of force generated
by particles placed at the zeroes and with masses proportional to the zeroes' multiplicities (see Section \ref{sec:location}). Later refinements are due to Jensen \cite{Je} and Lucas \cite{Lu}, but the first treatises on this matter are Marden's book \cite{Ma} and, primarily, Walsh's monograph \cite{Wa}
that collects most of the results on the number and location of critical points of complex polynomials and harmonic functions known at that date. In general dimension, even for harmonic functions, results are sporadic and rely on explicit formulae or symmetry arguments.
\par
Two well known questions in this context concern the location of the {\it hot spot} in a heat conductor --- a hot spot is a point of (absolute or relative) maximum temperature in the conductor. The situation described by \eqref{heat1}-\eqref{heat2} corresponds with the case of a {\it grounded} conductor. By some asymptotic analysis, under appropriate assumptions on $\fhi$, one can show that the hot spots {\it originate} from the set of maximum points of the function $d_\Om(x)$ --- the distance of $x\in\Om$ from $\Ga$ --- and tend to
the maximum points of the unique positive solution of \eqref{eigenfunction}, as $t\to\infty$.  In the case $\Om$ 
is convex, we have only one hot spot, as already observed. In Section \ref{sec:location}, we will describe three techniques to locate it; some of them extend their validity to locate the maximum points of the solutions to \eqref{torsion} and \eqref{eigenfunction}. We will also give an account of what it is known about convex conductors that admit a stationary hot spot (that is the hot spot does not move with time).
\par
It has also been considered the case in which the homogeneous Dirichlet boundary condition in \eqref{heat2} is replaced by the homogeneous Neumann condition:
\begin{equation}
\label{heat3}
u_\nu=0 \ \mbox{ on } \ \Ga\times(0,\infty).
\end{equation}
These settings describe the evolution of temperature in an {\it insulated} conductor of given constant 
initial temperature and has been made popular by a conjecture of J. Rauch \cite{Ra} that would imply that the hot spot must tend to a boundary point. Even if we now know that it is false for a general domain, the conjecture holds true for certain planar convex domains but it is still standing for unrestrained convex domains. 
\par
The remainder of the
paper is divided into three sections that reflect the aforementioned features. In Section \ref{sec:harmonic}, we shall describe the local properties of critical points of harmonic functions or, more generally, of solutions of elliptic equations, that lead to estimates of the size of critical sets.
In Section \ref{sec:number}, we shall focus on bounds for the number of critical points that depend on the boundary behavior of the relevant solutions and/or the geometry of $\Ga$. 
Finally, in Section \ref{sec:location}, we shall address the problem of locating the possible critical points.   As customary for a survey, our presentation will stress ideas rather than proofs. 
\par
This paper is dedicated with sincere gratitude to Giovanni Alessandrini ---  an inspiring mentor, a supportive colleague and  a genuine friend --- on the occasion of his $60^\mathrm{th}$ birthday. 
Much of the material presented here was either inspired by his ideas or actually carried out in his research with the author.

\section{The size of the critical set of a harmonic function}
\label{sec:harmonic}
A harmonic function in a domain $\Om$ is a solution of the Laplace's equation
$$
\De u=u_{x_1 x_1}+\cdots+u_{x_N x_N}=0 \ \mbox{ in } \ \Om.
$$
It is well known that harmonic functions are analytic, so there is no difficulty to define their critical points or the {\it critical set} 
$$
\cC(u)=\{x\in\Om: \na u(x)=0\}.
$$ 
Before getting into the heart of the matter, we present a relevant example.
\subsection{Harmonic polynomials}
In dimension two, we have a powerful tool since we know that a harmonic function is (locally) the real or imaginary part of a holomorphic function. 
This remark provides our imagination with a reach set of examples on which we can speculate. For instance, the harmonic function
$$
u=\re(z^n)=\re\left[(x+i y)^n\right], \ n\in\NN,
$$
already gives some insight on the properties of harmonic functions we are interested in.
In fact, we have that
$$
u_x-iu_y=n z^{n-1};
$$
thus, $u$ has only one distinct critical point, $z=0$, but it is more convenient to say that
$u$ has $n-1$ critical points at $z=0$ or that $z=0$ is a critical point with {\it multiplicity} $m$ with $m=n-1$. By virtue of this choice, we can give a topological meaning to $m$. 
\par
To see that, it is advantageous to represent $u$ in polar coordinates:
$$
u=r^n \cos(n\te);
$$
here, $r=|z|$ and $\te$ is the principal branch of $\arg z$, that is we are assuming that $-\pi\le\te<\pi$.  
Thus, the topological meaning of $m$ is manifest when we look at the level ``curve''
$\{ z: u(z)=u(0)\}$:  it is made of $m+1=n$ straight lines passing through the critical point $z=0$, divides the plane into $2n$ {\it cones} (angles), each of amplitude $\pi/n$
and the sign of $u$ changes across those lines (see Fig. 2.1). One can also show that the signed angle $\om$ formed by  $\na u$ and the direction of the positive real semi-axis, since it equals $-(n-1) \arg z$, increases by $2\pi m$ while $z$ makes a complete loop clockwise around $z=0$; thus, $m$ is a sort of {\it winding number} for $\na u$.


\begin{figure}[h]
\label{fig:power}
\centering
\includegraphics[scale=.9]{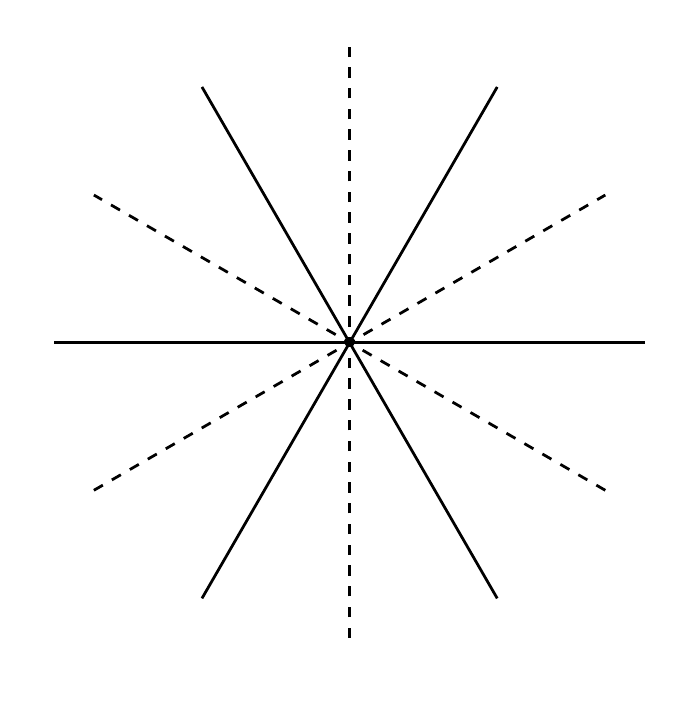}

\caption{Level set diagram of $u=r^6\cos(6\te)$ at the critical point $z=0$; $u$ changes sign from positive to negative at dashed lines and from negative to positive at solid lines.}
\end{figure}


The critical set of a homogeneous polynomial $P: \RR^N\to \RR$ is a cone in $\RR^N$. Moreover, if $P$ is also harmonic (and non-constant) one can show that 
\begin{equation}
\label{dimension}
\mbox{dimension of $\cC(u)$} \le N-2.
\end{equation}

\subsection{Harmonic functions}
\label{sub:harmonic}
If $N=2$ and $u$ is any harmonic function, the picture is similar to
that outlined in the example. In fact, we can again consider the ``complex gradient'' of $u$,
$$
g=u_x-iu_y,
$$
and observe that $g$ is {\it holomorphic} in $\Om$, since $\pa_{\ol{z}}g=0$, and hence analytic. Thus, the zeroes of $g$ (and hence the critical points of $u$) in $\Om$ are {\it isolated} and have
{\it finite multiplicity}. If $z_0$ is a zero with multiplicity $m$ of $g$, then we can write that
$$
g(z)=(z-z_0)^m h(z),
$$
where $h$ is holomorphic in $\Om$ and $h(z_0)\not=0$.
\par
On the other hand, we also know that $u$ is locally the real part of a holomorphic function $f$ 
and hence, since $f'=g$, by an obvious normalization, it is not difficult to infer that 
$$
f(z)=\frac1{n}\,(z-z_0)^{n}k(z),
$$
where $n=m+1$ and $k$ is holomorphic and $k(z_0)=h(z_0)\not=0$.
Passing to polar coordinates by $z=z_0+r e^{i\te}$ tells us that
$$
f(z_0+r e^{i\te})=\frac{|h(z_0)|}{n}\,r^{n} e^{i(n\te+\te_0)}+O(r^{n+1}) \ \mbox{ as } \ r\to 0,
$$
where $\te_0=\arg h(z_0)$.
Thus, we have that
$$
u=\frac{|h(z_0)|}{n}\,r^{n} \cos(n\te+\te_0)+O(r^{n+1}) \ \mbox{ as } \ r\to 0,
$$
and hence, modulo a rotation by the angle $\te_0$, in a small neighborhood of $z_0$, we can say that the critical level curve $\{z: u(z)=u(z_0)\}$ is very similar to that described in the example with $0$ replaced by $z_0$. In particular, it is made of $n$ simple curves passing through $z_0$ and any two adjacent curves meet 
at $z_0$ with an angle that equals $\pi/n$ (see Fig. 2.2).


\begin{figure}[h]
\label{fig:harmonic}
\centering
\includegraphics[scale=.9]{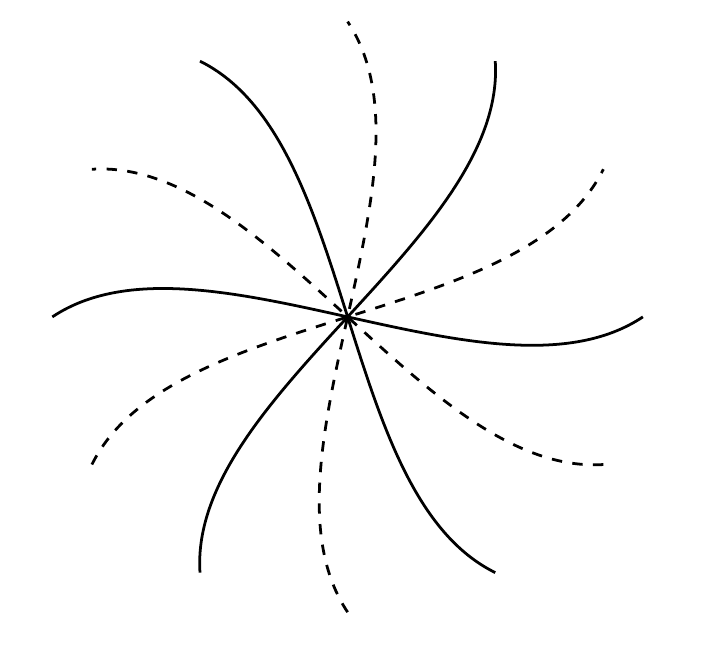}
\caption{Level set diagram of a harmonic function at a critical point with multiplicity $m=5$.
The curves meet with equal angles at the critical point.}
\end{figure}

If $N\ge 3$, similarly, a harmonic function can be approximated near 
a zero~$0$ by a homogeneous harmonic polynomial of some degree $n$:
\begin{equation}
\label{approximation}
u(x)= P_n(x)+O(|x|^{n+1}) \ \mbox{ as } \ |x|\to 0.
\end{equation}
However, the structure of the set $\cC(u)$ depends on whether $0$ is an isolated critical point of $P_n$ or not. In fact, if $0$ is not isolated, then $\cC(u)$ and $\cC(P_n)$ could not be diffeomorphic in general, as shown by the harmonic function
$$
u(x,y,z)=x^2-y^2+(x^2+y^2)\,z-\frac23\,z^3, \quad (x,y,z)\in\RR^3.
$$ 
Indeed, if $P_2(x,y,z)=x^2-y^2$, $\cC(P_2)$ is the $z$-axis, while $\cC(u)$ is made of $5$ isolated points (\cite{Pe}).
 
\subsection{Elliptic equations in the plane}
\label{subsec:elliptic}

These arguments can be repeated with some necessary modifications for solutions of uniformly elliptic
equations of the type \eqref{elliptic},
where the variable coefficients $a, b, c$ are Lipschitz continuous and $d, e$ are bounded measurable 
on $\Om$ and the uniform ellipticity is assumed to take the following form:
$$
ac-b^2=1 \ \mbox{ in } \ \Om.
$$
\par
Now, the classical theory of {\it quasi-conformal} mappings comes in our aid (see \cite{Be,Ve} and also \cite{AM1,AM2}). By the {\it uniformization theorem} (see \cite{Ve}), there exists a quasi-conformal mapping $\zi(z)=\xi(z)+i\,\eta(z)$, satisfying the
equation 
$$
\zi_{\ol{z}}=\kappa(z)\,\zi_z \ \mbox{ with } \ |\kappa(z)|=\frac{a+c-2}{a+c+2}<1,
$$ 
such that the function $U$ defined by $U(\zi)=u(z)$ satisfies the equation
$$
\De U+P\,U_\xi+Q\,U_\eta=0 \ \mbox{ in } \ \zi(\Om), 
$$
where $P$ and $Q$ are real-valued functions depending on the coefficients in \eqref{elliptic}  and are essentially bounded on $\zi(\Om)$. Notice that, since the composition of $\zi$ with a conformal mapping is still quasi-conformal, if it is convenient, by the Riemann mapping theorem, we can choose $\zi(\Om)$ to be the unit disk $\DD$.
\par
By setting $G=U_\xi-i\,U_\eta$, simple computations give that
$$
G_{\ol{\zi}}=R\,G+\ol{R}\,\ol{G} \ \mbox{ in } \DD,
$$
where $R=(P+i\,Q)/4$ is essentially bounded. This equation tells us that $G$ is 
a {\it pseudo-analytic} function for which the following {\it similarity principle} holds (see \cite{Ve}):
there exist two functions, $H(\zi)$ holomorphic in $\DD$ and $s(\zi)$  H\"older continuous on
the whole $\CC$, such that
\begin{equation}
\label{pseudo-analytic}
G(\zi)=e^{s(\zi)} H(\zi) \ \mbox{ for } \ \zi\in\DD.
\end{equation}
\par
Owing to \eqref{pseudo-analytic}, it is clear that the critical points of $u$, by means of the mapping $\zi(z)$, correspond to the 
zeroes of $G(\zi)$ or, which is the same, of $H(\zi)$ and hence we can claim that they are isolated 
and have a finite multiplicity.
\par
This analysis can be further extended if the coefficients $d$ and $e$ are zero, that is for the solutions of \eqref{elliptic2}. In this case, we can even assume
that the coefficients $a, b, c$ be merely essentially bounded on $\Om$, provided that we agree that
$u$ is a non-constant {\it weak} solution of \eqref{elliptic}. It is well known that, with these assumptions, solutions 
of \eqref{elliptic} are in general only H\"older continuous and the usual definition of critical point is
no longer possible. However, in \cite{AM2} we got around this difficulty by introducing a different notion of critical point, that is still consistent with the topological structure of the level curves of $u$ at its critical values.
\par
To see this, we look for a surrogate of the harmonic conjugate for $u$. In fact, \eqref{elliptic} implies
that the $1$-form
$$
\om=-(b\,u_x+c\,u_y)\,dx+(a\,u_x+b\,u_y)\,dy
$$
is closed (in the weak sense) in $\Om$ and hence, thanks to the theory developed in \cite{BN}, we can find a so-called {\it stream function} 
$v\in W^{1,2}(\Om)$ whose differential $dv$ equals $\om$, in analogy with the theory of gas dynamics (see \cite{BS}). 

Thus, in analogy with what we have done in Subsection \ref{sub:harmonic}, we find out that the function 
$f=u+i\,v$ satisfies the equation
\begin{equation}
\label{quasi-regular}
f_{\ol{z}}=\mu\,f_z
\end{equation}
where
$$
\mu=\frac{c-a-2ib}{2+a+c} \ \mbox{ and } \ |\mu|\le\frac{1-\la}{1+\la}<1 \ \mbox{ in } \ \Om,
$$
and $\la>0$ is a lower bound for the smaller eigenvalue of the matrix of the coefficients: 
$$
\left(
\begin{array}{cc}
a & b\\
b & c
\end{array}
\right).
$$
The fact that $f\in W^{1,2}(\Om,\CC)$ implies that $f$ is a {\it quasi-regular} mapping 
that can be factored as 
$$
f=F\circ\chi \ \mbox{ in } \ \Om,
$$
where $\chi:\Om\to\DD$ is a {\it quasi-conformal homeomorphism}  and $F$ is holomorphic in $\DD$
(see \cite{LV}). Therefore, the following {\it representation formula } holds:
$$
u=U(\chi(z)) \ \mbox{ for } \ z\in\Om,
$$
where $U$ is the real part of $F$.
\begin{figure}[h]
\label{fig:elliptic}
\centering
\includegraphics[scale=.8]{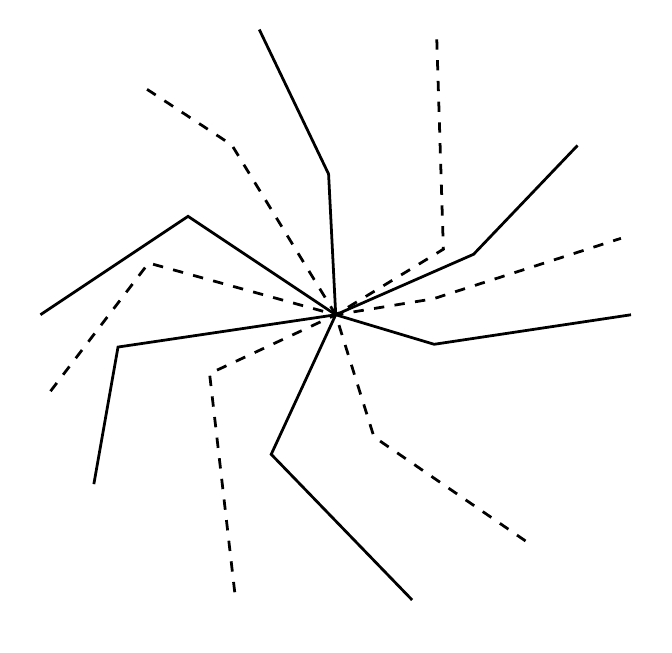}
\caption{Level set diagram of a solution of an elliptic equation with discontinuous coefficients at a geometric critical point with multiplicity $m=5$.  At that point, any two consecutive curves meet with positive angles, possibly not equal to one another.}
\end{figure}
\par
This formula informs us that the level curves of $u$ can possibly be distorted by the homeomorphism $\chi$, but preserve the topological structure of a harmonic function (see Fig. 2.3). This remark gives grounds to the 
definition introduced in \cite{AM2}: $z_0\in\Om$ is a {\it geometric critical point} of $u$ if the gradient of $U$ vanishes at $\chi(z_0)\in\DD$. In particular, geometric critical points are isolated and can be classified by a sort of multiplicity.

\subsection{Quasilinear elliptic equations in the plane}
A similar local analysis can be replicated when $N=2$ for quasilinear equations of type
$$
\dv\{A(|\na u|)\,\na u\}=0,
$$
where $A(s)>0$ and $0<\la\le 1+s\,A'(s)/A(s)\le\La$ for every $s>0$ and some constants $\la$ and $\La$.


\begin{figure}[h]
\label{fig:degenerate}
\centering
\includegraphics[scale=.85]{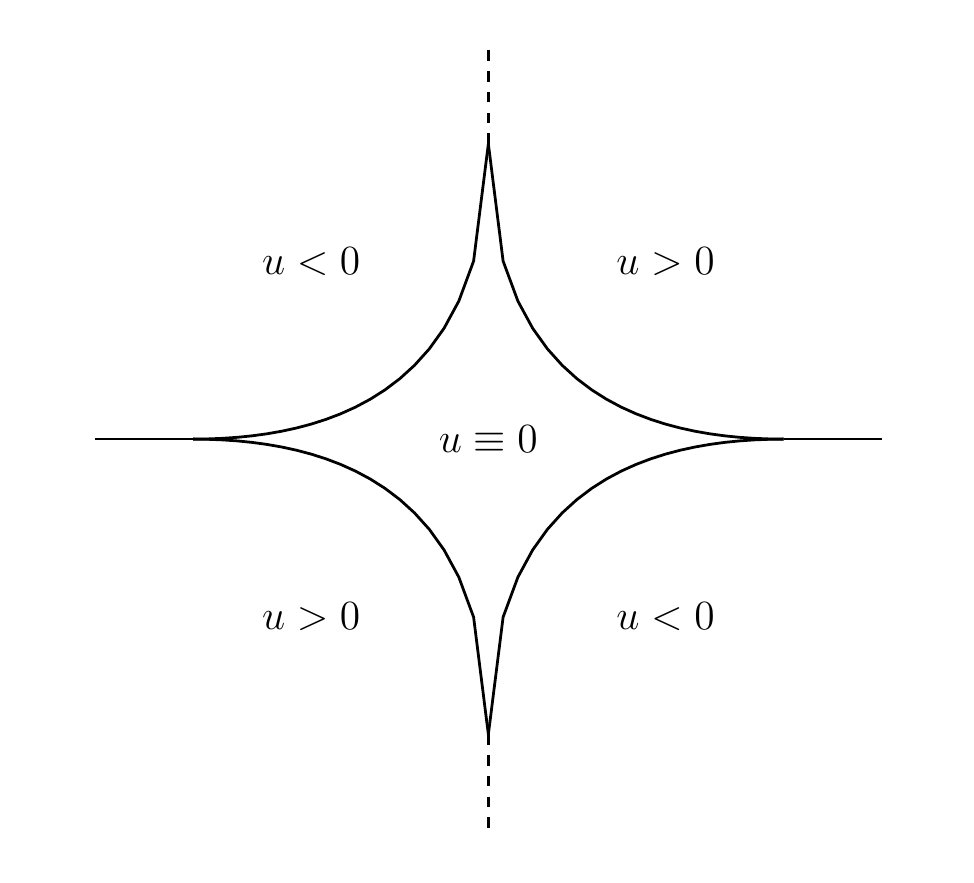}
\caption{Level set diagram of a solution of a degenerate quasilinear elliptic equation with $B(s)=\sqrt{1+s^2}$ at a critical value.}
\end{figure}

These equations can be even degenerate, such as the $p$-Laplace equation with $1<p<\infty$ (see \cite{AR}). It is worth mentioning that also the case in which $A(s)=B(s)/s$, where $B$ is increasing, with $B(0)>0$, and superlinear and growing polynomially at infinity (e.g. $B(s)=\sqrt{1+s^2}$), has been studied in~\cite{CM}. In this case the function $1+s\,A'(s)/A(s)$ vanishes at $s=0$ and it turns out that the critical points of a solution $u$ (if any) are {\it never} isolated (Fig. 2.4).

\subsection{The case $\bf N\ge 3$}
\label{sub:circles}
As already observed, critical points of harmonic functions in dimension $N\ge 3$ may not be isolated. 
Besides the example given in Section \ref{sub:harmonic}, another concrete example is given by the function
$$
u(x,y,z)=J_0\Bigl(\sqrt{x^2+y^2}\Bigr)\,\cosh(z), \quad (x,y,z)\in\RR^3,
$$
where $J_0$ is the first Bessel function: the gradient of $u$ vanishes at the origin and on the circles on the plane $z=0$ having radii equal to the zeroes of the second Bessel function $J_1$. It is clear that
a region $\Om$ can be found such that $\cC(u)\cap\Om$ is a {\it bounded continuum}.

\smallskip

Nevertheless, it can be proved that $\cC(u)$ always has locally finite $(N-2)$-dimensional Hausdorff measure $\cH^{N-2}$. A nice argument to see this was suggested to me by D. Peralta-Salas \cite{Pe}.
If $u$ is a non-constant harmonic function and we suppose that $\cC(u)$ has dimension $N-1$,  then the general theory of analytic sets implies that there is an open and dense subset of $\cC(u)$ which is an analytic sub-manifold (see \cite{KP}). Since $u$ is constant on a connected component of the
critical set, it is constant on $\cC(u)$, and its gradient vanishes. Thus,  by the
Cauchy-Kowalewski theorem $u$ must be constant in a neighborhood of $\cC(u)$, and hence everywhere by unique continuation. Of course, this argument would also work for solutions of
an elliptic equation of type
\begin{equation}
\label{elliptic-N}
\sum_{i, j=1}^N a_{ij}(x)\,u_{x_i x_j}+\sum_{j=1}^N b_j(x)\,u_{x_j}=0 \ \mbox{ in } \ \Om,
\end{equation}
with analytic coefficients.
\par
When the coefficients $a_{ij}, b_j$ in \eqref{elliptic-N} are of class $C^\infty(\Om)$, the result has been proved in \cite{HHHN} (see also \cite{Ha}): if $u$ is a non-constant solution of \eqref{elliptic-N}, then for any compact subset $K$
of $\Om$ it holds that
\begin{equation}
\label{hausdorff}
\cH^{N-2}(\cC(u)\cap K)<\infty.
\end{equation}
The proof is based on an estimate similar to \eqref{dimension} for the complex
dimension of the singular set in $\CC^N$ of the complexification of the polynomial $P_n$ in the approximation \eqref{approximation}. 
\par
The same result does not hold for solutions of equation
\begin{equation}
\label{elliptic-complete}
\sum_{i, j=1}^N a_{ij}(x)\,u_{x_i x_j}+\sum_{j=1}^N b_j(x)\,u_{x_j}+c(x)\,u=0 \ \mbox{ in } \ \Om,
\end{equation}
with $c\in C^\infty(\Om)$. For instance the gradient of the first Laplace-Dirichlet
eigenfunction for a spherical annulus vanishes exactly on a $(N-1)$-dimensional sphere. A more general counterexample is the following
(see \cite[Remark p.~362]{HHHN}): let $v$ be of class $C^\infty$ and with non-vanishing gradient in the unit ball $B$ in $\RR^N$; the function $u=1+v^2$ satisfies the equation
$$
\De u-c u=0 \ \mbox{ with } \ c =\frac{\De v^2}{1+v^2}\in C^\infty(B);
$$
we have that $\cC(u)=\{ x\in B: v(x)=0\}$ and it has been proved that any closed subset of $\RR^N$ can
be the zero set of a function of class $C^\infty$ (see \cite{To}). 
\par
However, once \eqref{hausdorff} is settled, it is rather easy to show that the {\it singular set} 
$$
\cS(u)=\cC(u)\cap u^{-1}(0)=\{x\in\Om: u(x)=0, \na u(x)=0\}
$$
of a non-constant solution of \eqref{elliptic-complete}
also has locally finite $(N-2)$-dimensional Hausdorff measure \cite[Corollary 1.1]{HHHN}. 
This can be done by a trick, since around any point in $\Om$ there always exists a {\it positive} solution $u_0$ of \eqref{elliptic-complete} and it turns out that  the function $w=u/u_0$ is a solution of an equation like \eqref{elliptic-N} and that $\cS(u)\subseteq\cC(w)$.
In particular the set of critical points on the nodal line of an eigenfunction of
the Laplace operator has locally finite $(N-2)$-dimensional Hausdorff measure.
\par
Nevertheless, for a solution of \eqref{elliptic-N} the set $\cS(u)$ can be very complicated, as a simple example in \cite[p. 361]{HHHN}) shows: the function $u(x,y,z)=xy+f(z)^2$, where
$f$ is a smooth function with $|f f''|+(f')^2<1/4$ that vanishes exactly on an arbitrary given closed subset $K$ of $\RR$, is a solution of
$$
u_{xx}+u_{yy}+u_{zz}-(f^2)''(z)\,u_{xy}=0 \ \mbox{ and } \ \cS(u)=\{(0,0)\}\times K.
$$
\par
Heuristically, as in the $2$-dimensional case, the proof of \eqref{hausdorff} is essentially based on the observation that, by Taylor's expansion, a harmonic function $u$ can be approximated near any of its zeroes by a homogeneous harmonic polynomial $P_m(x_1,\dots, x_n)$ of degree $m\ge 1$. Technically, the authors use  the fact that the complex dimension of the critical set in $\CC^N$ of the complexified polynomial $P_m(z_1,\dots, z_N)$ is bounded by $N-2$. A $C^\infty$-perturbation argument and an inequality from geometric measure theory then 
yield that, near a zero of $u$, the $\cH^{N-2}$-measure of $\cC(u)$ can be bounded in terms of $N$ and $m$. The extension of these arguments to the case of a solution of 
\eqref{elliptic-N} is then straightforward. 
Recently in \cite{CNV}, \eqref{hausdorff} has been extended to the case of solutions
of elliptic equations of type
$$
\sum_{i, j=1}^N \{a_{ij}(x)\,u_{x}\}_{x_j}+\sum_{j=1}^N b_j(x)\,u_{x_j},
$$ 
where the coefficients $a_{ij}(x)$ and $b_j(x)$ are assumed to be Lipschitz continuous and essentially bounded, respectively.

\section{The number of critical points}
\label{sec:number}

A more detailed description of the critical set $\cC(u)$ of a harmonic function $u$ can be obtained if
we assume to have some information on its behavior on the boundary $\Ga$ of $\Om$. 
While in Section \ref{sec:harmonic} the focus was on a qualitative description of the set $\cC(u)$, here we are concerned with establishing bounds on the number of critical points. 

\subsection{Counting the critical points of a harmonic function in the plane}
An exact counting formula is given by the following result.

\begin{thm}[\cite{AM1}]
\label{th:am1}
Let $\Om$ be a bounded domain in the plane and let 
$$
\Ga=\bigcup_{j=1}^J \Ga_j,
$$
where $\Ga_j, j=1,\dots, J$ are simple closed curves of class $C^{1,\al}$. Consider a
harmonic function $u\in C^1(\ol{\Om})\cap C^2(\Om)$ that satisfies the Dirichlet boundary condition
\begin{equation}
\label{condenser}
u=a_j \ \mbox{ on } \  \Ga_j, \,j=1,\dots, J,
\end{equation} 
where $a_1, \dots, a_J$ are given real numbers, not all equal.
\par
Then $u$ has in $\ol{\Om}$ a finite number of critical points $z_1,\dots, z_K$; if 
$m(z_1), \dots,$ $m(z_K)$ denote their multiplicities, then the following identity holds:
\begin{equation}
\label{identity}
\sum_{z_k\in\Om}m(z_k)+\frac12\,\sum_{z_k\in\Ga}m(z_k)=J-2.
\end{equation} 
 \end{thm}
 \begin{figure}[h]
\label{fig:capacitor}
\centering
\includegraphics[scale=.75]{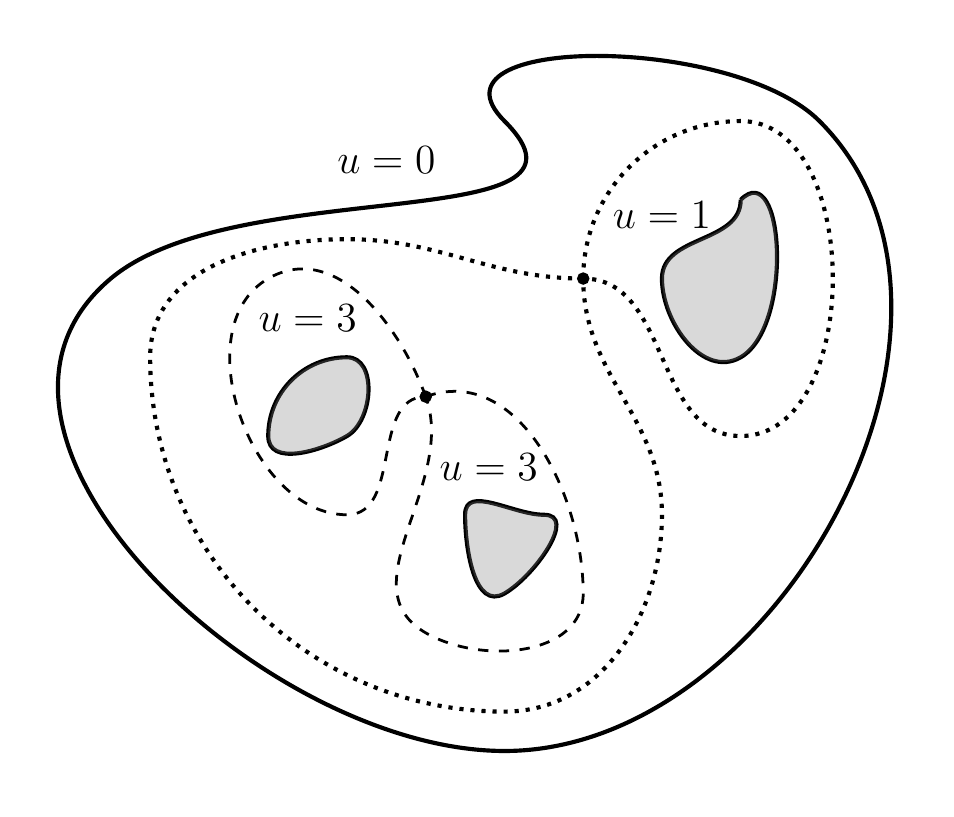}
\caption{An illustration of Theorem \ref{th:am1}: the domain $\Om$ has $3$ holes; $u$ has exactly $2$ critical points; dashed and dotted are the level curves at critical values.}
\end{figure}
\par
Thanks to the analysis presented in Subsection \ref{subsec:elliptic}, this theorem
still holds if we replace the Laplace equation in \eqref{condenser} by
 the general elliptic equation \eqref{elliptic}. In fact, modulo a suitable change of
 variables, we can use \eqref{pseudo-analytic} with $\im(s)=0$ on the boundary.
\par
The function considered in Theorem \ref{th:am1} can be interpreted in physical terms as the potential
in an electrical capacitor and hence its critical points are the points of equilibrium
of the electrical field (Fig. 3.1). 
\par
 The proof of Theorem \ref{th:am1} relies on the fact that the critical points of $u$ are the zeroes of the
 holomorphic function $f=u_x-i\,u_y$ and hence they can be counted with their multiplicities by applying the classical {\it argument principle} to $f$ with some necessary modifications. The important remark is that, since the boundary components are level curves for $u$, the gradient of $u$ is parallel on them to the (exterior) unit normal $\nu$ to the boundary, and hence $\arg f=-\arg\nu$. 
\par
Thus, the situation is clear if $u$ does not have critical points on $\Ga$: the argument principle gives at once that
\begin{multline*}
 \sum_{z_k\in\Om}m(z_k)=\frac1{2\pi i}\int_{+\Ga}\frac{f'(z)}{f(z)}\,dz\\
=\frac1{2\pi}\,\inc(\arg f,+\Ga)=
\frac1{2\pi}\,\inc(-\arg \nu,+\Ga)=
-[1- (J-1)]=J-2,
 \end{multline*}
 where by $\inc(\cdot,+\ga)$ we intend the {\it increment} of an angle on an oriented curve $+\ga$ and by $+\Ga$ we mean that $\Ga$ is trodden in such a way that  $\Om$
 is on the left-hand side.
\par
 If $\Ga$ contains critical points, we must first prove that they are also isolated.
 This is done, by observing that, if $z_0$ is a critical point belonging to some component $\Ga_j$, since $u$ is constant on $\Ga_j$, by the {\it Schwarz's reflection principle} (modulo a conformal transformation of $\Om$), $u$ can be extended to 
 a function~$\widetilde{u}$ which is harmonic in a whole neighborhood of $z_0$. Thus,
 $z_0$ is a zero of the holomorphic function $\widetilde{f}=\widetilde{u}_x-i\,\widetilde{u}_y$ 
 and hence is isolated and with finite multiplicity. Moreover, the increment of $\arg\widetilde{f}$ on an oriented closed simple curve~$+\ga$ around $z_0$ is exactly twice as much as that of
 $\arg f$ on the part of $+\ga$ inside $\Om$. This explains the second addendum in \eqref{identity}.
 
 \medskip
 
%

Notice that condition \eqref{condenser} can be re-written as
$$
u_\tau=0 \ \mbox{ on } \ \Ga,
$$
where $\tau:\Ga\to\SS^1$ is the {\it tangential} unit vector field on $\Ga$.
We cannot hope to obtain an identity as \eqref{identity} if $u_\tau$ is not constant. 
However, a bound for the number of critical points of a harmonic function (or a solution of \eqref{elliptic}) can be derived in a quite general setting. 
\par
In what follows,
we assume that $\Om$ is as in Theorem \ref{th:am1} and that
$\ell:\Ga\to\SS^1$ denotes a (unitary) vector field of class $C^1(\Ga, \SS^1)$ of
given topological degree~$D$, that can be defined as
\begin{equation}
\label{degree}
2\pi\,D=\inc(\arg(\ell),+\Ga).
\end{equation}
 Also, we will use the following definitions:
\begin{enumerate}[(i)]
\item
if $(\cJ^+,\cJ^-)$ is a decomposition of $\Ga$ into two disjoint subsets such that $u_\ell\ge 0$ on $\cJ^+$ and $u_\ell\le 0$ on $\cJ^-$, we denote by 
$M(\cJ^+)$ the number of connected components of $\cJ^+$ which are {\it proper subsets} of some component $\Ga_j$ of $\Ga$ and set:
$$
M=\min\{ M(\cJ^+): (\cJ^+, \cJ^+) \mbox{ decomposes } \Ga \};
$$
\item
if $\cI^\pm=\{z\in\Ga: \pm\,u_\ell(z)>0\}$, by $M^\pm$ we denote the number
of connected components of $\cI^\pm$ which are {\it proper subsets} of some component $\Ga_j$ of $\Ga$.
\end{enumerate}
Notice that in (i) the definition of $M$ does not change if we replace $\cJ^+$ by $\cJ^-$.

\begin{thm}[\cite{AM1}]
\label{th:am1b}
Let $u\in C^1(\ol{\Om})\cap C^2(\Om)$ be harmonic in $\Om$ and
denote by $m(z_j)$ the multiplicity of a zero $z_j$ of $f=u_x-i\,u_y$.
\begin{enumerate}[(a)]
\item
If $M$ is finite and $u$ has no critical point in $\Ga$, then
$$
\sum_{z_j\in\Om}m(z_j)\le M-D;
$$
\item
if $M^++M^-$ is finite, then 
$$
\sum_{z_j\in\Om}m(z_j)\le \left[\frac{M^++M^-}{2}\right]-D,
$$
where $[x]$ is the greatest integer $\le x$.
\end{enumerate}
\end{thm}

This theorem is clearly less sharp than Theorem \ref{th:am1} since, in that setting, it does not give information about critical points on the boundary. However, it gives the same information on the number of interior critical points, since in the setting of Theorem \ref{th:am1} the degree of the field $\tau$ on $+\Ga$ equals $2-J$ and $M=0$. 


\begin{figure}[h]
\label{fig:oblique}
\centering
\includegraphics[scale=.7]{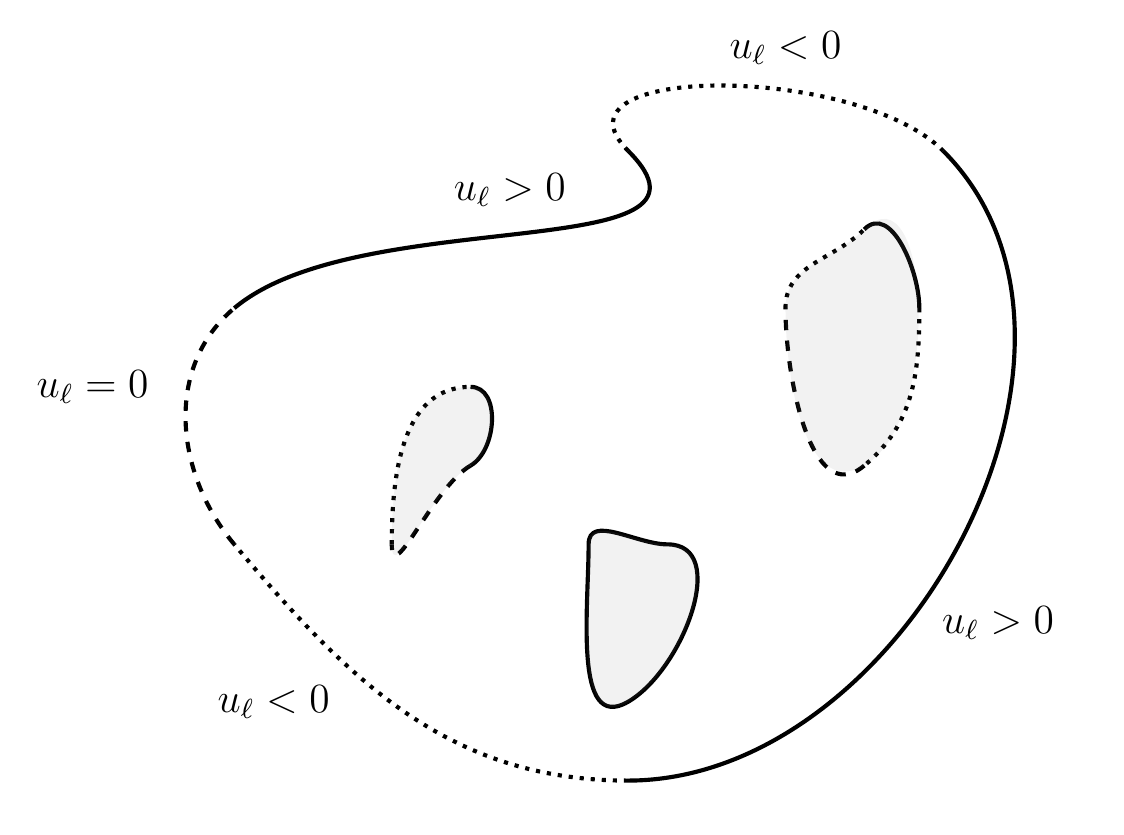}
%
%
%
%
%
\caption{An illustration of Theorem \ref{th:am1b}. Here, $M=M^+=4$; $M^-=4$; $D=-2$
if $\ell=\nu$ or $\tau$; $D=1$ if $\ell=z/|z|$ and the origin is in $\Om$; $D=0$ if $\ell=(1,0)$ or $(0,1)$. }
\end{figure}


\par
The possibility of choosing the vector field $\ell$ arbitrarily makes Theorem \ref{th:am1b} a very flexible tool: for instance, the number of critical points in $\Om$ can be estimated from information on the tangential, normal, co-normal, partial, or radial (with respect to some origin) derivatives (see Fig. 3.2).  
\par
As an illustration, it says that in a domain topologically equivalent to a disk, in order to have $n$ interior critical point the normal (or tangential, or co-normal) derivative of a harmonic function must change sign at least $n+1$ times and a partial derivative at least $n$ times. Thus, Theorem \ref{th:am1b} helps 
to choose Neumann data that insures the absence of critical points in $\Om$. For this reason,  in its general form for elliptic operators, it has been useful in the study of EIT and other similar inverse problems.

\medskip

We give a sketch of the proof of (a) of Theorem \ref{th:am1b}, that hinges on the simple fact that, if we set $\te=\arg(\ell)$ and 
$\om=\arg(u_x-i\,u_y)$, then
$$
u_\ell=\ell\cdot\na u=|\na u|\,\cos(\te+\om).
$$
Hence, if $(\cJ^+, \cJ^-)$ is a minimizing decomposition of $\Ga$ as in (i), then
$$
|\om+\te|\le\frac{\pi}{2} \ \mbox{ on } \  \cJ^+ \ \mbox{ and } \ |\om+\te-\pi|\le\frac{\pi}{2} \ \mbox{ on } \  \cJ^-.
$$ 
\par
Thus, two occurrences must be checked. If a component $\Ga_j$ is contained in $\cJ^+$ or $\cJ_-$,
then 
$$
\left|\frac1{2\pi}\,\inc(\om+\te,+\Ga_j)\right|\le\frac12,
$$
that implies that $\om$ and $-\te$ must have the same increment, being the right-hand side an integer.  If $\Ga_j$ contains points of both $\cJ^+$ and $\cJ^-$, instead, if $\si^+\subset\cJ^+$ and
$\si^-\subset\cJ^-$ are two consecutive components on $\Ga_j$, then 
$$
\frac1{2\pi}\,\inc(\om+\te,+(\si^+\cap\si^-)\le 1.
$$
Therefore, if $M_j$ is the number of connected components  of $\cJ^+\cap\Ga_j$ (which equals that of $\cJ^+\cap\Ga_j$), then 
$$
\frac1{2\pi}\,\inc(\om+\te,+\Ga_j)\le M_j,
$$
and hence
\begin{multline*}
\sum_{z_k\in\Om}m(z_k)=\frac1{2\pi}\,\inc(\om,+\Ga)=
\frac1{2\pi}\,\inc(\om+\te,+\Ga)-D\\
=\sum_{j=1}^J\frac1{2\pi}\,\inc(\om+\te,+\Ga_j)\le 
\sum_{j=1}^J M_j-D=M-D.
\end{multline*}

\bigskip

{\bf The obstacle problem.} An estimate similar to that of Theorem \ref{th:am1b} has been obtained also for $N=2$ by Sakaguchi \cite{Sa}
for the {\it obstacle problem}. Let $\Om$ be bounded and simply connected and let $\psi$ be a given function in $C^2(\ol{\Om})$ --- the obstacle. There exists a unique solution $u\in H^1_0(\Om)$ such that $u\ge\psi$ in $\Om$ of the obstacle problem 
$$
\int_\Om \na u\cdot\na(v-u)\,dx\ge 0 \ \mbox{ for every } \  v\in H^1_0(\Om) \ \mbox{ such that } \ u\ge\psi.
$$ 
It turns out that $u\in C^{1,1}(\ol{\Om})$ and $u$ is harmonic outside of the {\it contact set}
$I=\{x\in\Om: u(x)=\psi(x)\}$. 
In \cite{Sa} it is proved that, if the number of connected components of local maximum points of $\psi$ equals $J$, then
$$
\sum_{z_k\in\Om\setminus I}m(z_k)\le J-1,
$$
with the usual meaning for $z_k$ and $m(z_k)$. In \cite{Sa}, this result is also shown to hold for
a more general class of quasi-linear equations. The proof of this result is based on the 
analysis of the level sets of $u$ at critical values, in the wake of~\cite{Al1} and \cite{HW}.

\bigskip

Topological bounds as in Theorems \ref{th:am1} or \ref{th:am1b} are not possible in dimension greater than $2$. We give two examples.


\begin{figure}[h]
\label{fig:torus}
\centering
\includegraphics[scale=.7]{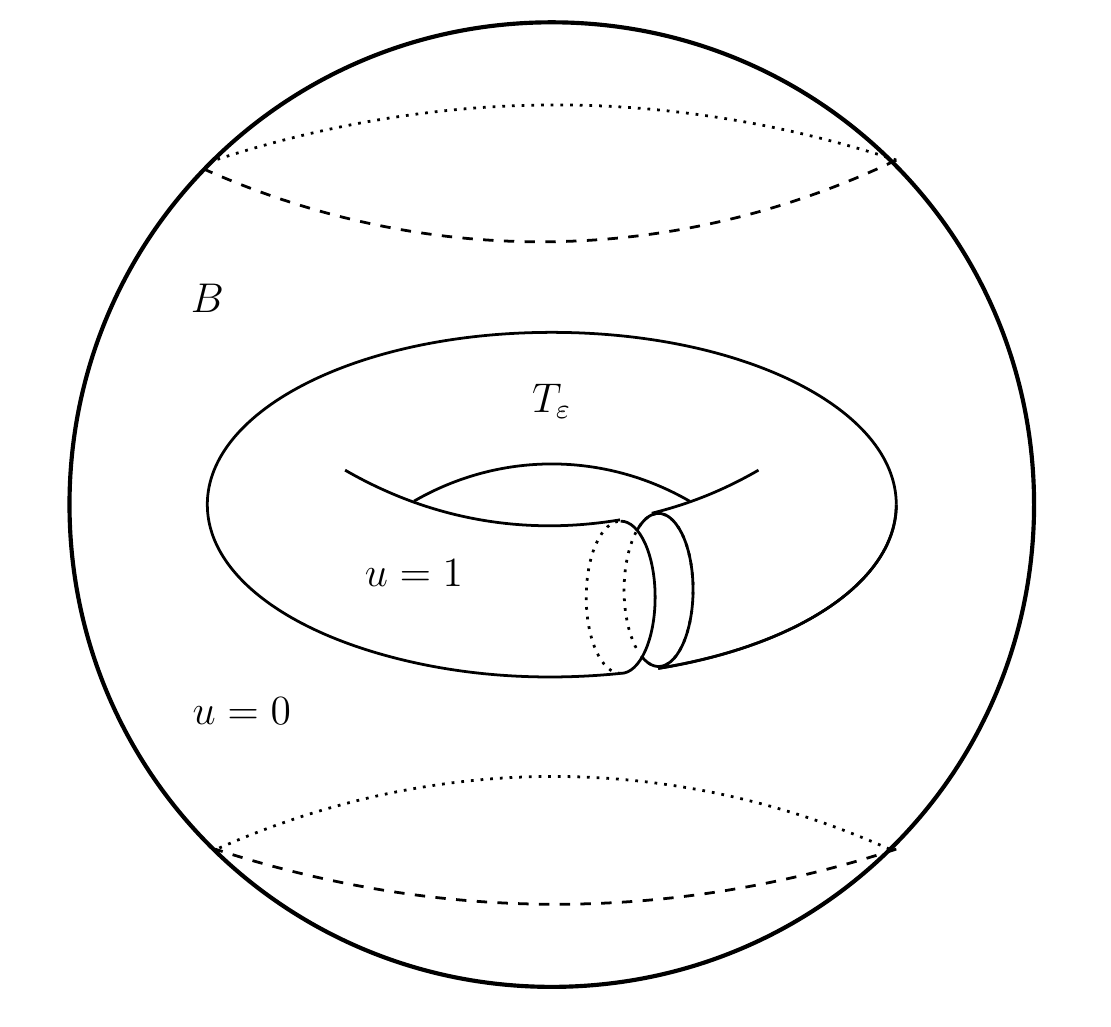}
%
%
%
\caption{The broken doughnut in a ball: $u$ must have a critical point near the center of $B$ 
and one between the ends of $T$.}
\end{figure}


{\bf The broken doughnut in a ball.}
The first is an adaptation of one contained in \cite{EP1} and reproduces the situation of 
Theorem \ref{th:am1} (see Fig. 3.3). Let $B$ be the unit ball centered at the origin in $\RR^3$ and $T$ an open torus  with center of symmetry at the origin and such that $\ol{T}\subset B$. We can always choose 
 coordinate axes in such a way that the $x_3$-axis is the axis of revolution for $T$ 
 and hence define the set $T_\ve=\{ x\in T: x_2<\ve^{-1}|x_1|\}$.
 $\ol{T_\ve}$ is simply connected and tends to $T$ as $\ve\to 0^+$. 
Now, set $\Om_\ve=B\setminus\ol{T_\ve}$ and 
consider a capacity potential for $\Om$, that is the harmonic function in $\Om_\ve$ with the following boundary values
$$
u=0 \ \mbox{ on } \ \pa B, \quad u=1 \ \mbox{ on } \pa T_\ve.
$$
Since $\Om_\ve$ has $2$ planes of symmetry (the $x_1x_2$ and $x_2x_3$ planes), the partial derivatives $u_{x_1}$ and $u_{x_3}$ must be zero 
on the two segments that are the intersection of $\Om_\ve$ with the $x_2$-axis. If
 $\si$ is the segment that contains the origin, the restriction of $u$ to $\ol{\si}$ equals
 $1$ at the point $\ol{\si}\cap\pa T_\ve$, is $0$ at the point $\ol{\si}\cap\pa B$, is bounded at the origin by a constant $<1$ independent of $\ve$, and can be made arbitrarily close to $1$ between the ``ends'' of $T_\ve$, when $\ve\to 0^+$,  It follows that,
 if $\ve$ is sufficiently small,  $u_{x_2}$ (and hence $\na u$) must
 vanish twice on $\si$. 
 \par
 It is clear that this argument does not depend on the size or on small deformations of $T$. 
Thus, we can construct in $B$ a (simply connected) ``chain''~$C_\ve$ of an arbitrary number $n$ of such tori,  by gluing them together: the solution in the domain obtained by replacing $T_\ve$ by $C_\ve$ will then have at least $2n$ critical points. 

\medskip

{\bf Circles of critical points.}
The second example shows that, in general dimension, a finite number of sign changes of some derivative of a harmonic function $u$ on the boundary does not even imply that $u$ has
a finite number of critical points.
\par 
To see this, consider the harmonic function is Subsection \ref{sub:circles}:
$$
u(x,y,z)=J_0(\sqrt{x^2+y^2})\,\cosh(z).
$$
It is easy to see that, for instance, on any sphere centered at the origin the normal derivative $u_\nu$
changes its sign a finite number of times. However, if the radius of the sphere is larger than the first positive 
zero of $J_1=0$, the corresponding ball contains at least one circle of critical points.

\medskip

{\bf Star-shaped annuli.}  Nevertheless, if some additional geometric information is added, something can be done. Suppose that 
$\Om=D_0\setminus\ol{D_1}$, where $D_0$ and $D_1$ are two domains in $\RR^N$, with
boundaries of class $C^1$ and such that 
$\ol{D_1}\subset D_0$. Suppose that $D_0$ and $D_1$ 
are  {\it star-shaped} with respect to the same origin $O$ placed in $D_1$, that is the segment $OP$
is contained in the domain for every point $P$ chosen in it. Then, the {\it capacity potential} $u$
defined as the solution of the Dirichlet problem
$$
\De u=0 \ \mbox{ in } \ \Om, \quad u=0 \ \mbox{ on } \ \pa D_0, \quad u=1 \ \mbox{ on } \ \pa D_1,
$$
does not have critical points in $\ol{\Om}$. This is easily proved by considering the harmonic function 
$$
w(x)=x\cdot\na u(x), \ x\in\Om.
$$ 
Since $D_0$ and $D_1$ are starshaped and of class $C^1$, $w\ge 0$ on $\pa\Om$. By the strong maximum principle, then $w>0$ in $\Om$; in particular, $\na u$ does not vanish in~$\Om$ and all the sets $D_1\cup \{x\in\ol{\Om}: u(x)>s\}$ turn out to be star shaped too (see \cite{Ev}). This theorem can be extended to the capacity potential defined in $\Om=\RR^N\setminus\ol{D_1}$ as the solution of
$$
\De u=0 \ \mbox{ in } \ \Om, \quad u=1 \ \mbox{ on } \ \pa\Om, \quad u\to 0 \ \mbox{ as } \ |x|\to\infty.
$$
Such results have been extended in \cite{Fc,Pu,Sl1} to a very general class of nonlinear elliptic equations.

\subsection{Counting the critical points of Green's functions on manifolds}

With suitable restrictions on the coefficients, \eqref{elliptic2} can be regarded as the {\it Laplace-Beltrami equation} on the Riemannian surface $\RR^2$ equipped with the metric
$$
c\,(dx)^2-2 b\,(dx)(dy)+a\,(dy)^2.
$$
Theorems \ref{th:am1} and \ref{th:am1b} can then be interpreted accordingly.
\par
This point of view has been considered in a more general context in \cite{EP2,EP3},
where the focus is on Green's functions of a $2$-dimensional complete Riemannian surface $(M, g)$ of finite topological type (that is, the first fundamental group of $M$ is finitely generated). A Green's function is a symmetric 
function $\cG(x,y)$ that satisfies in $M$ the equation
\begin{equation}
\label{laplace-beltrami}
-\De_g \cG(\cdot,y)=\de_y ,
\end{equation}
where $\De_g$ is the Laplace-Beltrami operator induced by the metric $g$ and $\de_y$ is the {\it Dirac delta} centered at a point $y\in M$. 
\par
A symmetric Green's function $\cG$ can always be constructed 
by an approximation argument introduced in \cite{LT}:  an increasing sequence of compact subsets~$\Om_n$ containing $y$ and exhausting $M$ is introduced and $\cG$ is then defined as the limit on compact subsets of $M\setminus\{ y\}$ of the sequence $\cG_n-a_n$, where $\cG_n$ is the solution of 
\eqref{laplace-beltrami} such that $\cG_n=0$ on $\Ga_n$ and $a_n$ is a suitable constant.
A Green's function defined in this way is generally not unique, but has many properties in common with the fundamental solution for Laplace's equation in  the Euclidean plane. 
\par
With these premises, in \cite{EP2,EP3} it has been proved the following notable topological bound:
$$
\mbox{number of critical points of $\cG$ } \le 2\mathfrak{g}+\mathfrak{e}-1,
$$
where $\mathfrak{g}$ and $\mathfrak{e}$ are the {\it genus} and the {\it number of ends} of $M$;
the number $2\mathfrak{g}+\mathfrak{e}-1$ is known as the first {\it Betti number} of $M$. Moreover, if the Betti number is attained, then $\cG$ is {\it Morse}, that is at its critical points the Hessian matrix is non-degenerate. In \cite{EP2}, it is also shown that, in dimensions greater than two, an upper bound by topological invariants is impossible.
\par
Two different proofs are constructed in \cite{EP2} and \cite{EP3}, respectively. Both proofs 
are based on the following {\it uniformization principle}: since $(M, g)$ is a smooth manifold of
finite topological type, it is well known (see \cite{KY}) that there exists a compact surface $\Si$
endowed with a metric $g'$ of constant curvature, a finite number $J\ge 0$ of isolated points 
points $p_j\in\Si$ and a finite number $K\ge 0$ of (analytic) topological disks $D_k\subset\Si$ such
that $(M, g)$ is conformally isometric to the manifold $(M', g')$, where $M'$ is interior of 
$$
M'=\Si\setminus\left(\bigcup_{j=1}^J\{ p_j\} \cup \bigcup_{k=1}^K D_k\right).
$$
That means that there exist a diffeomorphism $\Phi: M\to M'$ and a positive function $f$ on $M$ such that $\Phi^* g'=f g$; it turns out that the genus $\mathfrak{g}$ of $\Si$ and the number $J+K$ --- that equals the number $\mathfrak{e}$ ends of $M$ ---
determine $M$ up to diffeomorphisms.

\medskip

The proof in \cite{EP2} then proceeds by analyzing the transformed Green's function $\cG'=\cG\circ\Phi^{-1}$. It is proved that $\cG'$ satisfies the problem
$$
-\De_{g'} \cG'(\cdot, y')=\de_{y'}-\sum\limits_{j=1}^J c_j\de_{p_j}\ \mbox{ in the interior of } \  M', \qquad
\cG'=0 \mbox{ on } \ \bigcup_{k=1}^K \pa D_k,
$$
where $y'=\Phi(y)$ and the constants $c_j$, possibly zero (in which case $\cG'$ would be $g'$-harmonic near $p_j$), sum up to $1$. Thus,
a local blow up analysis of the {\it Hopf index} $\mathfrak{I}(z_n), j=1,\dots, N$, of the gradient of $\cG'$ at the critical points $z_1,\dots, z_N$ (isolated and with finite multiplicity), together
with the {\it Hopf Index Theorem} (\cite{Mo,MC}), yield the formula
$$
\sum_{n=1}^N\mathfrak{I}(z_n)+\sum_{c_j\not=0}\mathfrak{I}(p_j)=\chi(\Si^*),
$$
where $\chi(\Si^*)$ is the Euler characterstic of the manifold
$$
\Si^*=\Si\setminus\left( D_{y'}\cup\bigcup_{k=1}^K D_k \right)
$$
and $D_{y'}$ is a sufficiently small disk around $y'$. Since $\chi(\Si^*)$ is readily computed as
$1-2\mathfrak{g}-K$ and $\mathfrak{I}(z_n)\le -1$, one then obtains that
\begin{multline*}
\mbox{ number of critical points of  $\cG'$}=-\sum_{n=1}^N\mathfrak{I}(z_n)=\\
2\mathfrak{g}+K-1+\sum_{c_j\not=0}\mathfrak{I}(p_j)\le 2\mathfrak{g}+J+K-1=2\mathfrak{g}+\mathfrak{e}-1.
\end{multline*}
Of course, the gradient of $\cG^*$ vanishes if and only if that of $\cG$ does.
 
\medskip
 
The proof contained in \cite{EP3} has a more geometrical flavor and focuses on the study of the 
integral curves of the gradient of $\cG$. This point of view is motivated by the fact that in Euclidean space the Green's function (the fundamental solution) arises as the electric potential of a charged particle at $y$, so that its critical points correspond to equilibria and the integral curves of its gradient field are the lines of force classically studied in the XIX century. Such a description relies on techniques of dynamical systems rather than on the toolkit of partial differential equations.
\par
We shall not get into the details of this proof, but we just mention that it gives a more satisfactory 
portrait of the integral curves connecting the various critical points of $\cG$ --- an issue that has rarely been studied.
 
\subsection{Counting the critical points of eigenfunctions}

The bounds and identities on the critical points that we considered so far are based on 
a crucial topological tool: the index $\mathfrak{I}(z_0)$ of a critical point $z_0$. 
\par
For a function
$u\in C^1(\Om)$, the integer 
$\mathfrak{I}(z_0)$ is the {\it winding number or degree} of the vector field $\na u$ around $z_0$ and is related to the portrait of the set $\cN_u=\{z\in\cU: u(z)=u(z_0)\}$ for a sufficiently small neighborhood $\cU$ of $z_0$. As a matter of fact, 
if $z_0$ is an isolated critical point of $u$, one can distinguish two situations (see \cite{AM1,Ro}):
\begin{enumerate}[(I)]
\item 
if $\cU$ is sufficiently small, $\cN_u=\{z_0\}$ and $\mathfrak{I}(z_0)=1$;
\item
if $\cU$ is sufficiently small, $\cN_u$ consists of $n$ simple curves and, if $n\ge 2$, each pair of such curves  crosses at $z_0$ only; it turns out that $\mathfrak{I}(z_0)=1-n$.
\end{enumerate}
Critical points with index $\mathfrak{I}$ equal to $1$, $0$, or negative are called 
{\it extremal, trivial}, or {\it saddle} points, respectively (see \cite{AM1}) . A saddle point is {\it simple} or {\it Morse} if the hessian matrix of $u$ at that point is not trivial.
\par
In the cases we examined so far, we always have that $\mathfrak{I}(z_0)\le -1$, that is $z_0$ is a saddle point,  since (I) and (II) with $n=1$ cannot
occur, by the maximum principle.
\par
The situation considerably changes when $u$ is a solution of \eqref{torsion}, \eqref{eigenfunction}, or \eqref{heat1}. Here, we shall give an account of what can be said for solutions of \eqref{eigenfunction}. The same ideas can be used for solutions of the semilinear equation
$$
-\De u=f(u) \ \mbox{ in } \ \Om,
$$
subject to a homogeneous Dirichlet boundary condition, where the non-linearity $f:\RR\to\RR$
satisfies the assumptions:
$$
f(t)>0 \ \mbox{ if } t>0 \quad\mbox{ or }\quad f(t)/t>0 \ \mbox{ for } \ t\not=0
$$
(see \cite{AM1} for details).
We present here the following result that is in the spirit of Theorem \ref{th:am1}. 

\begin{thm}[\cite{AM1}]
\label{th:am1c}
Let $\Om$ be as in Theorem \ref{th:am1} and $u\in C^1(\ol{\Om})\cap C^2(\Om)$ be a solution of \eqref{eigenfunction}. If $z_0\in\ol{\Om}$ is an isolated critical point of $u$ in $\ol{\Om}$, then
\begin{enumerate}[(A)]
\item
either $z_0$ is a nodal critical point, that is $z_0\in\cS(u)$, and the function $u_x-i\,u_y$ is asymptotic
to $c\,(z-z_0)^m$, as $z\to z_0$, for some $c\in\CC\setminus\{ 0\}$ and $m\in\NN$,
\item
or $z_0$ is an extremal, trivial, or simple saddle critical point. 
\end{enumerate}
\par
Finally, if all the critical points of $u$ in $\ol{\Om}$ are isolated\,\footnote{This assumption can be removed when $\Om$ is simply connected, by using the analyticity of $u$ (see \cite{AM1})}, the following identity holds:
\begin{equation}
\label{identity2}
\sum_{z_k\in\Om}m(z_k)+\frac12\,\sum_{z_k\in\Ga}m(z_k)+n_S-n_E=J-2.
\end{equation}
Here, $n_S$ and $n_E$ denote the number of the simple saddle and extremal points of $u$.  
\end{thm}
Thus, a bound on the number of critical points in topological terms is not possible --- additional information of different nature should be added.

\medskip

The proof of this theorem can be outlined as follows. 
\par
First, one observes that, at a nodal critical point
$z_0\in\Om$, $\De u$ vanishes, and hence the situation described in Subsection \ref{sub:harmonic} is in order,
that is $u_x-i\,u_y$ actually behaves as specified in (A) and the index $\mathfrak{I}(z_0)$ equals $-m$. If $z_0\in\Ga$, a reflection argument like the one used for Theorem \ref{th:am1} can be used, so that $z_0$ can be treated as an interior nodal critical point of an extended function with vanishing laplacian at $z_0$ and (A) holds; in this case, however, as done for Theorem \ref{th:am1}, the
contribution of $z_0$ must be counted as $-m/2$.
\par
Secondly, one examines non-nodal critical points. At these points $\De u$ is either positive or negative. If, say, $\De u(z_0)<0$, then at least one eigenvalue of the hessian matrix of $u$ must be negative and
the remaining eigenvalue is either positive (and hence a simple saddle point arises), negative (and hence a maximum point arises) or zero (and hence, with a little more effort, either a trivial or a simple saddle point arises). Thus, the total index of these points sums up to $n_E-n_S$.
\par
Finally, identity \eqref{identity2} is obtained by applying Hopf's index theorem in a suitable manner.

\subsection{Extra assumptions: the emergence of geometry}

As emerged in the previous subsection, topology is not enough to control the number of critical points
of an eigenfunction or a torsion function. Here, we will explain how some geometrical information
about $\Om$ can be helpful.
\par
Convexity is a useful information. If the domain $\Om\subset\RR^N$, $N\ge 2$, is convex, one can expect that the solution $\tau$ of
\eqref{torsion} and the only {\it positive} solution $\phi_1$ of \eqref{eigenfunction} --- it exists and, as is well known, corresponds to the first Dirichlet eigenvalue $\la_1$ --- have only one critical point (the maximum point). This expectation is realistic, but a rigorous proof is not straightforward. 
\par
In fact, one has to first show that $\tau$ and $\phi_1$ are
{\it quasi-concave}, that is one shows the convexity of the level sets 
$$
\{ x\in\Om: u(x)\ge s\} \ \mbox{ for every } \ 0\le s\le \max_{\ol{\Om}} u,
$$
for $u=\tau$ or $u=\phi_1$.
It should be noted that $\phi_1$ is {\it never concave} and
examples of convex domains $\Om$ can be constructed such that $\tau$ is not concave (see \cite{Ko}). 
\par 
The quasi-concavity of $\tau$ and $\phi_1$ can be proved in several different ways (see 
\cite{BiS,BL,CS,IS,Ke,Ko,Sl2}). Here, we present the argument used in \cite{Ko}. There, the desired quasi-convexity is obtained by showing that the functions $\si=\sqrt{\tau}$ and
$\psi=\log\phi_1$ are concave functions ($\tau$ and
$\phi_1$ are then said {\it $1/2$-concave} and {\it  log-concave}, respectively). 
\par
In fact, one shows that  $\si$ and $\psi$ satisfy the conditions
$$
\De\si=-\frac{1+2\,|\na\si|^2}{2\si} \ \mbox{ in } \ \Om, \quad \si=0 \ \mbox{ on } \Ga,
$$
and
$$
\De\psi=-(\la_1+|\na\psi|^2) \ \mbox{ in } \ \Om, \quad \psi=-\infty \ \mbox{ on } \Ga.
$$
The concavity test established by Korevaar in \cite{Ko}, based on a maximum principle for
the so-called {\it concavity function} (see also \cite{Ka}), applies to these two problems and guarantees that both
$\si$ and $\psi$ are concave. With similar arguments, one can also prove that the solution
of \eqref{heat1}-\eqref{heat2} is $\log$-concave in $x$ for any fixed time $t$.
\par
The obtained quasi-concavity implies in particular that, for $u=\tau$ or $\phi_1$, the set of critical points $\cC(u)$, that here coincides
with the set
$$
\cM(u)=\Bigl\{ x\in\Om: u(x)=\max_{\ol{\Om}} u\Bigr\},
$$
is convex. 
This set cannot contain more than one point, due to the analyticity of $u$. In fact, if it
contained a segment, being the restriction of $u$ analytic on the chord of $\ol{\Om}$ containing that segment, 
$u$ would be a {\it positive} constant on this chord and this is impossible, since $u=0$ at the endpoints of this chord.
\par
This same argument makes sure that, if $\fhi\equiv 1$ in a convex domain $\Om$, then for any fixed $t>0$ there is 
a {\it unique} point $x(t)\in\Om$ --- the so-called {\it hot spot} --- at which the solution of \eqref{heat1}-\eqref{heat2} attains its maximum in $\ol{\Om}$, that is 
$$
h(x(t),t)=\max_{x\in\ol{\Om}} h(x,t) \ \mbox{ for } \ t>0.
$$
The location of $x(t)$ in $\Om$ will be one of the issues in the next section.

\medskip

{\bf A conjecture.} Counting (or estimating the number of) the critical points of $\tau$, $\phi_1$, or $h$ when $\Om$ is not convex seems a difficult task.
For instance, to the author's knowledge, it is not even known whether or not the uniqueness of the maximum
point holds true if $\Om$ is assumed to be {\it star-shaped} with respect to some origin.
\par
We conclude this subsection by offering and justifying a conjecture on the number of hot spots 
in a bounded simply connected domain $\Om$ in $\RR^2$. 
To this aim, we define for $t>0$ the set of hot spots as
$$
\cH(t)=\{ x\in\Om: x \mbox{ is a local maximum point of $h(\cdot,t)$}\}.
$$
\par
We shall suppose that the function $\fhi$ in \eqref{heat2} is continuous, non-negative and not identically equal to zero in $\Om$, so that, by  {\it Hopf's boundary point lemma}, $\cH(t)\cap\Ga=\varnothing$. Also, by an argument 
based on the analyticity of $h$ similar to that 
used for the uniqueness of the maximum point in a convex domain, we can be sure that $\cH(t)$ is made of isolated points (see \cite{AM1} for details). (A parabolic version of ) Theorem \ref{th:am1c} then yields that
$$
n_E(t)-n_S(t)=1,
$$
where $n_E(t)$ and $n_S(t)$ are the number extremal and simple saddle points of $h(\cdot,t)$;  clearly $n_E(t)$ is the cardinality of $\cH(t)$. An estimate on the total number of critical points
of $h(\cdot,t)$ will then follow from one on $n_E(t)$.
\par
Notice that, if $\la_n$ and $\phi_n$, $n\in\NN$, are Dirichlet eigenvalues (arranged in increasing order) and eigenfunctions (normalized in $L^2(\Om)$) of the Laplace's operator in $\Om$, then the following {\it spectral formula}
\begin{equation}
\label{spectral}
h(x,t)=\sum_{n=1}^\infty \widehat{\fhi}(n)\,  \phi_n(x) e^{-\la_n t} \ \mbox{ holds for } \ x\in\ol{\Om} \ \mbox{ and } \ t>0,
\end{equation}
where $\widehat{\fhi}(n)$ is the Fourier coefficient of $\fhi$ corresponding to $\phi_n$. 
Then we can infer that 
$e^{\la_1 t} h(x,t)\to\widehat{\fhi}(1)\,\phi_1(x)$ as $t\to\infty$, with
$$
\widehat{\fhi}(1)=\int_\Om \fhi(x)\,\phi_1(x)\,dx>0,
$$
and the convergence is uniform on $\ol{\Om}$ under sufficient assumptions on $\fhi$ and $\Om$.
This information implies that, if $x(t)\in\cH(t)$, then
\begin{equation}
\label{large times}
\dist(x(t),\cH_\infty)\to 0 \ \mbox{ as } \ t\to\infty,
\end{equation}
where $\cH_\infty$ is the set of local maximum points of $\phi_1$. 
\par
Now, our conjecture concerns the influence of the shape of $\Om$ on the number $n_E(t)$. To rule out the possible influence of the values of $\fhi$, we assume that $\fhi\equiv 1$: then we know that there holds the following asymptotic formula (see~\cite{Va}):
\begin{equation}
\label{varadhan}
\lim_{t\to 0^+} 4t\,\log[1-h(x,t)]=-d_\Ga(x)^2 \ \mbox{ for } \ x\in\ol{\Om};
\end{equation}
here, $d_\Ga(x)$ is the {\it distance} of a point $x\in\ol{\Om}$ from the boundary $\Ga$. The convergence in \eqref{varadhan} is uniform on $\ol{\Om}$ under suitable regularity assumptions on $\Ga$.

\begin{figure}[h]
\label{fig:cocoon}
\centering
\includegraphics[scale=.8]{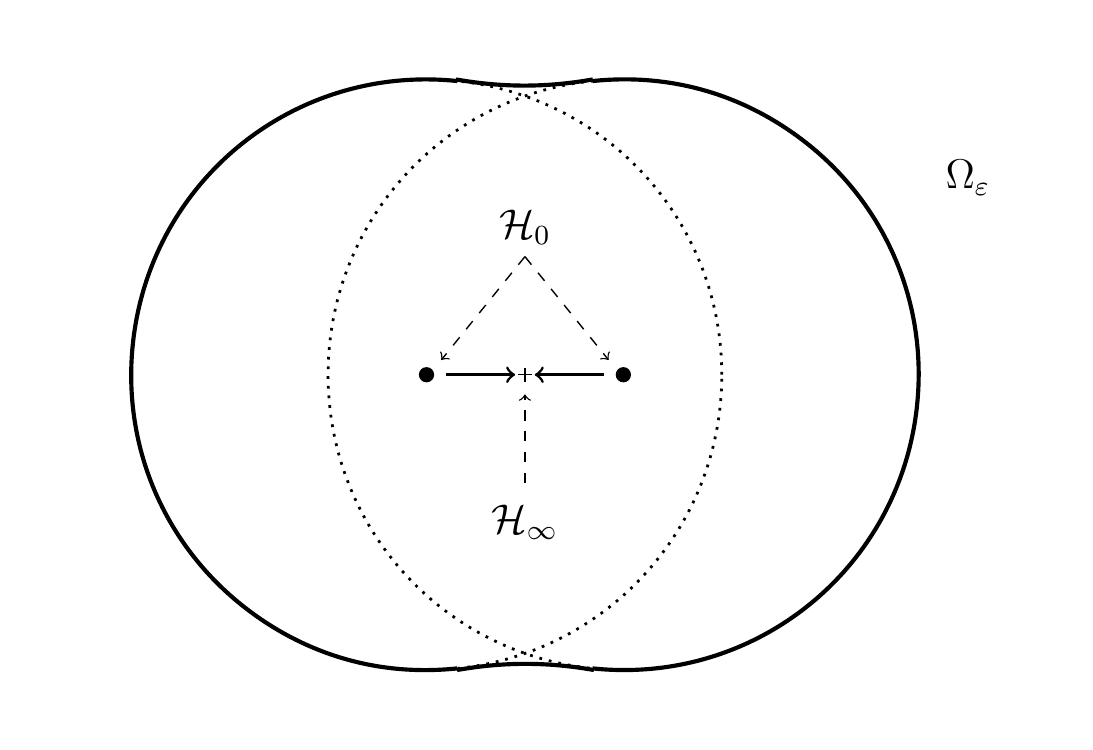}
\caption{As time $t$ increases, $\cH(t)$ goes from $\cH_0$, the set of maximum points of $d_\Ga$,
to $\cH_\infty$, the set of maximum points of $\phi_1$.}
\end{figure}

\par
Now, suppose that $d_\Ga$ has {\it exactly} $m$ distinct local (strict) maximum points in $\Om$. 
Formula \eqref{varadhan} suggests that, when $t$ is sufficiently small, $h(\cdot,t)$ has the same number $m$ of  maximum points in $\Om$. As time $t$ increases, one expects that the maximum points of $h(\cdot,t)$ do not increase in number. Therefore, the following bounds should hold:
\begin{equation}
\label{guess}
n_E(t)\le m \ \mbox{ and hence } \ n_E(t)+n_S(t)\le 2m-1 \ \mbox{ for every } \ t>0.
\end{equation}
 From the asymptotic analysis performed on \eqref{spectral}, we also derive that the {\it total
 number of critical points of} $\phi_1$ does not exceeds $2m-1$.
\par
We stress that \eqref{guess} cannot always hold with the equality sign. In fact, if $D_\ve^\pm$ denotes the unit disk centered at $(\pm\ve,0)$ and we consider the domain
$\Om_\ve$ obtained from $D_\ve^+\cup D_\ve^-$ by ``smoothing out the corners'' (see Fig. 3.4), we
notice that $m=2$ for every $0<\ve<1$, while $\Om_\ve$ tends to the unit ball centered at the origin and hence,  if $\ve$ is small enough, $\phi_1$ has only one critical point, being $\Om_\ve$ ``almost convex''.   
\par
Based on a similar argument, inequalities like \eqref{guess} should also hold for the number of critical points of the torsion function $\tau$. In fact, if $U_s$ is the solution of the one-parameter family of problems
$$
-\De U_s+s\,U_s=1 \ \mbox{ in } \ \Om, \quad U_s=0 \mbox{ on } \ \Ga,
$$
where $s$ is a positive parameter, we have that
$$
\lim_{s\to 0^+} U_s=\tau \quad \mbox{and} \quad \lim_{s\to\infty} \frac1{\sqrt{s}}\,\log[1-s\,U_s]=-d_\Ga,
$$ 
uniformly on $\ol{\Om}$ (see again \cite{Va}).
\par
We finally point out that the asymptotic formulas presented here hold in any dimension; thus, the bounds in  \eqref{guess} may be generalized in some way. 

\subsection{A conjecture by S. T. Yau}
To conclude this section about the number of critical points of solutions of partial differential equations, we cannot help mentioning a conjecture proposed in \cite{Ya} (also see \cite{EJN,JaN,JNT}). This is motivated
by the study of eigenfunctions of the Laplace-Beltrami operator $\De_g$ in a compact Riemannian
manifold $(M,g)$. 
\par
Let $\{\phi_k\}_{k\in\NN}$ be a sequence of eigenfunctions,
$$
\De_g\phi_k+\la_k\phi_k=0 \ \mbox{ in } \ M.
$$ 
Let $x_k\in M$ be a point of maximum for $\phi_k$ in $M$ and $B_k$ a geodesic ball centered at $x_k$ and with radius $C/\sqrt{\la_k}$. If we blow up $B_k$ to the unit disk in $\RR^2$ and let
$u_k/\max\phi_k$ be the eigenfunction after that change of variables, then a subsequence of $\{ u_k\}_{k\in\NN}$ will converge to a solution $u$ of
\begin{equation}
\label{yau}
\De u+u=0, \ |u|<1  \ \mbox{ in } \ \RR^2.
\end{equation}
If we can prove that $u$ has {\it infinitely many} isolated critical points, then we can expect that their number be unbounded also for the sequence $\{\phi_k\}_{k\in\NN}$.
\par
A naive insight built up upon the available concrete examples of entire eigenfunctions 
(the {\it separated} eigenfunctions in rectangular or polar coordinates) may suggest that it would be
enough to prove that  any solution of \eqref{yau} has infinitely many nodal domains. It turns out that
this is not always true, as a clever counterexample  obtained in \cite[Theorem 3.2]{EJN} shows:
{\it there exists a solution of \eqref{yau} with exactly two nodal domains}. 
\par
The counterexample is constructed by perturbing the solution of \eqref{yau} 
$$
f=J_1(r)\sin\te, 
$$
where $(r,\te)$ are the usual polar coordinates and $J_1$ is the second Bessel's function; $f$ has infinitely many nodal domains. The desired example is thus obtained by the perturbation $h=f+\ve\,g$,
where $g(x,y)=f(x-\de_x,y-\de_y)$ and $(\de_x,\de_y)$ is suitably chosen.  As a result, if $\ve$
is sufficiently small, the set $\{ (x,y)\in\RR^2: h(x,y)\not=0\}$ is made of two interlocked spiral-like domains (see \cite[Figure 3.1]{EJN}).
\par
A related result was proved in \cite{EP4}, where it is shown that there is no topological upper bound for the number of critical points of the first eigenfunction on Riemannian manifolds (possibly with
boundary) of dimension larger than two. In fact, with no restriction on the topology of the manifold, it is possible to construct metrics whose first
eigenfunction has as many isolated critical points as one wishes. 
\par
Recently, it has been proved  in \cite{JZ1} that, if $(M, g)$ is a non-positively curved surface with concave boundary, the number of nodal domains of $\phi_k$ diverges along a subsequence of eigenvalues of density $1$  (see also \cite{JZ2} for related results). The surface needs not have any symmetries. The number can also be shown to grow like $\log \la_k$ (\cite{Ze}). In light of such results,
Yau's conjecture was updated as follows: show that, for any (generic)
$(M,g)$ there exists at least one sub-sequence of eigenfunctions for
which the number of nodal domains (and hence of the critical points) tends to infinity (\cite{Ya2,Ze}).

\section{The location of critical points}
\label{sec:location}

\subsection{A little history}
The first result that studies the critical points of a function is probably {\it Rolle's theorem}:
between two zeroes of a differentiable real-valued function there is {\it at least one} critical point. 
Thus, a function that has $n$ distinct zeroes also has at least $n-1$ critical points --- an estimate from below --- and we roughly know where they are located.
\par
After Rolle's theorem, the first general result concerning the zeroes of the derivative of a general
polynomial is {\it Gauss's theorem}: if 
$$
P(z)=a_n\,(z-z_1)^{m_1}\cdots ,(z-z_K)^{m_K}, \ \mbox{ with } \ m_1+\cdots+m_K=n,
$$
is a polynomial of degree $n$, then
$$
\frac{P'(z)}{P(z)}=\frac{m_1}{z-z_1}+\cdots+\frac{m_K}{z-z_K}
$$
and hence the zeroes of $P'(z)$ are, in addition to the multiple zeroes of $P(z)$ themselves, the
roots of
$$
\frac{m_1}{z-z_1}+\cdots+\frac{m_K}{z-z_K}=0.
$$
These roots can be interpreted as the {\it equilibrium points} of the gravitational field generated
by the masses $m_1, \dots, m_K$ placed at the points $z_1, \dots, z_K$, respectively.

\begin{figure}[h]
\label{fig:lucas}
\centering
\includegraphics[scale=.9]{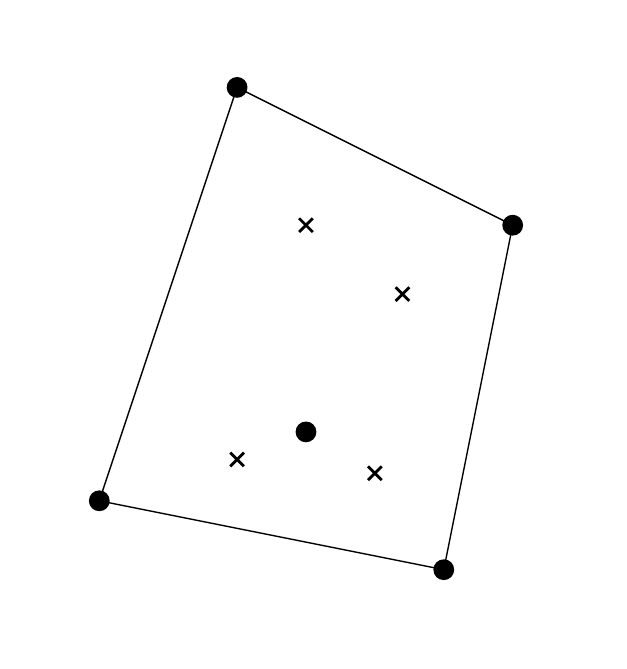}
\caption{Lucas's theorem: the zeroes of $P'(z)$ must fall in the convex envelope of those of $P(z)$.}
\end{figure}

\par
If the zeroes of $P(z)$ are placed on the real line then, by Rolle's theorem,  it is not difficult
to convince oneself that the zeroes of $P'(z)$ lie in the smallest interval of the real axis that contains
the zeroes of $P(z)$. This simple result has a geometrically expressive generalization in 
{\it Lucas's theorem}: the zeroes of $P'(z)$ lie in the {\it convex hull} $\Pi$ of the set $\{ z_1,\dots, z_K\}$ --- named {\it Lucas's polygon} ---and no such zero lies on $\pa\Pi$ unless is a multiple  zero $z_k$ of $P(z)$ or all the zeroes of $P(z)$ are collinear (see Fig. 4.1).
\par
 In fact, it is enough to observe that, if $z\notin\Pi$ or $z\in\pa\Pi$,
then all the $z_k$ lie in the closed half-plane $H$ containing them and the side of $\Pi$ which is the closest to $z$. Thus, if $\ell=\ell_x+i\,\ell_y$ is an outward direction to $\pa H$, we have that
$$
\re\left[\left(\sum_{k=1}^K \frac{m_k}{z-z_k}\right)\ell\right]=
\sum_{k=1}^K m_k\,\frac{\re\bigl[\ol{(z-z_k)}\,\ell\bigr]}{|z-z_k|^2}>0,
$$
since all the addenda are non-negative and not all equal to zero,
unless the $z_k$'s are collinear.

If $P(z)$ has real coefficients, we know that its non-real zeroes occur in conjugate pairs. Using
the circle whose diameter is the segment joining such a pair --- this is called a {\it Jensen's circle} of $P(z)$ ---
one can obtain a sharper estimate of the location of the zeroes of $P'(z)$: each non-real zero of $P'(z)$ lies on or within a Jensen's circle of $P(z)$. This result goes under the name of {Jensen's theorem} (see \cite{Wa} for a proof).
\par
All these results can be found in Walsh's treatise \cite{Wa}, that contains many other results about 
zeroes of complex polynomials or rational functions and their extensions to critical points of harmonic functions: among them restricted versions of Theorem \ref{th:am1} give information (i) on the critical points of the Green's function of an infinite region delimited by a finite collection of simple closed curves and (ii) of  harmonic measures generated by collections of Jordan arcs. Besides the {\it argument's principle}
already presented in these notes,  a useful ingredient used in those extensions is a {\it Hurwitz's theorem} (based on the classical {\it Rouch\'e's theorem}): if $f_n(z)$ and $f(z)$ are holomorphic in a
domain $\Om$, continuous on $\ol{\Om}$, $f(z)$ is non-zero on $\Ga$ and $f_n(z)$ converges
uniformly to $f(z)$ on $\ol{\Om}$, then there is a $n_0\in\NN$ such that, for $n>n_0$,  $f_n(z)$
and $f(z)$ have the same number of zeroes in $\Om$.

\subsection{Location of critical points of harmonic functions in space}

The following result is somewhat an analog of Lucas's theorem and is related to \cite[Theorem 1, p. 249]{Wa}, which holds in the plane.

\begin{thm}[\cite{EP1}] 
\label{ep1}
Let $D_1, \dots, D_J $ be bounded domains in $\RR^N$, $N\ge 3$, with boundaries of class $C^{1,\al}$ and
with mutually disjoint closures, and set
$$
\Om=\RR^N\setminus\bigcup_{j=1}^J \ol{D_j}.
$$
\par
Let $u\in C^0(\ol{\Om})\cap C^2(\Om)$ be the solution of the boundary value problem
\begin{equation}
\label{capacity}
\De u=0 \ \mbox{ in } \ \Om, \quad u=1 \ \mbox{ on } \ \Ga, \quad u(x)\to 0 \ \mbox{ as } |x|\to\infty.
\end{equation}
\par
If $\cK$ denotes the convex hull of   
$$
\bigcup_{j=1}^J D_j,
$$
then $u$ does not have critical points in $\ol{\RR^N\setminus\cK}$ (sse Fig. 4.2). 
\end{thm}

This theorem admits at least two proofs and it is worth to present  both of them.
The former is somewhat reminiscent of Lucas's proof and is based on an explicit formula for $u$,
$$
u(x)=\frac1{(N-2)\, \om_N}\int_\Ga \frac{u_\nu(y)}{|x-y|^{N-2}}\,dS_y, \ x\in\Om,
$$
that can be derived as a consequence of {\it Stokes's formula}. Here, $\om_N$ is the
surface area of a unit sphere in $\RR^N$, $dS_y$ denotes the $(N-1)$-dimensional surface measure, and $u_\nu$ is the (outward) normal derivative of $u$. 
\par
By the Hopf's boundary point lemma, $u_\nu>0$ on $\Ga$. Also, if $x\in\ol{\RR^N\setminus\cK}$, we can choose a hyperplane $\pi$
passing through $x$ and supporting $\cK$ (at some point). If $\ell$ is the unit
vector orthogonal to $\pi$ at $x$ and pointing into the half-space containing $\cK$, we have that
$(x-y)\cdot\ell$ is non-negative and is not identically zero for $y\in\Ga$. Therefore,
$$
u_\ell(x)=-\frac1{\om_N}\int_\Ga \frac{u_\nu(y)\,(x-y)\cdot\ell}{|x-y|^{N}}\,dS_y<0,
$$
which means that $\na u(x)\not=0$.


\begin{figure}[h]
\label{fig:ep}
\centering
\includegraphics[scale=.9]{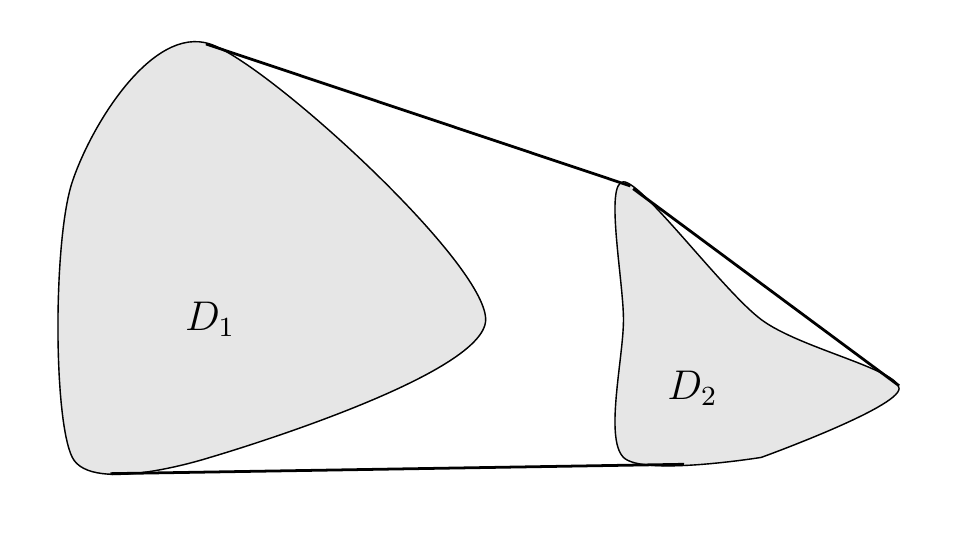}
\caption{No critical points outside of the convex envelope.}
\end{figure}


The latter proof is based on a symmetry argument (\cite{Sa2}) and, as it will be clear, can also be extended to more general non-linear equations. Let $\pi$ be any hyperplane contained in $\ol{\Om}$ and let $H$ be the open half-space containing $\cK$ and such that $\pa H=\pi$. Let $x'$ be the mirror reflection in $\pi$ of any point $x\in H\cap\Om$. Then the function defined by
$$
u'(x)=u(x') \ \mbox { for } \ x\in H\cap\Om
$$
is harmonic in $H\cap\Om$, tends to $0$ as $|x|\to\infty$ and 
$$
u'<u \ \mbox{ in } \ H\cap\Om, \quad u'=u \ \mbox{ on } \ \pi\setminus\Ga.
$$
Therefore, by the Hopf's boundary point lemma, $u_\ell(x)\not=0$ at any $x\in\pi\setminus\Ga$ for any direction $\ell$ not parallel to $\pi$. Of course, if $x\in\Ga\cap\pi$, we obtain that $u_\nu(x)>0$ by directly using the Hopf's boundary point lemma.

\medskip

Generalizations of Lucas's theorem hold for other problems. Here, we mention the well known result of Chavel and Karp \cite{CK} for the minimal solution of the Cauchy problem for the heat equation in
a Riemannian manifold $(M,g)$:
\begin{equation}
\label{ck}
u_t=\De_g u \ \mbox{ in } \ M\times(0,\infty), \qquad u=\fhi \ \mbox{ on } \ M\times\{ 0\},
\end{equation}
where  $\fhi$ is a bounded initial data with compact support in $M$. In \cite{CM}, it is shown that, if $M$ is complete, simply connected and of constant curvature, then the set of the {\it hot spots} of $u$,
$$
\cH(t)=\left\{x\in M: u(x,t)=\max_{y\in M}u(y,t)\right\},
$$
is contained in the convex hull of the support of $\fhi$. The proof is based on an explicit 
formula for $u$ in terms of the initial values $\fhi$. For instance, when $M=\RR^N$, we have the
formula
$$
u(x,t)=(4\pi t)^{-N/2}\int_{\RR^N} e^{-|x-y|^2}\fhi(y)\,dy \ \mbox{ for } \ (x,t)\in\RR^N\times(0,\infty).
$$
With this formula in hand, by looking at the second derivatives of $u$, one can also prove that there
is a time $T>0$ such that, for $t>T$,  $\cH(t)$ reduces to the single point
$$
\frac{\int_{\RR^N} y\,\fhi(y) dy}{\int_{\RR^N} \fhi(y) dy},
$$
which is the center of mass of the measure space $(\RR^N, \fhi(y) dy)$ (see \cite{JS}).
\par
We also mention here the work of Ishige and Kabeya (\cite{IK1,IK2,IK3}) on the large time behavior of  hot spots for solutions of the heat equation with a rapidly decaying potential and for the 
Schr\"odinger equation.

\subsection{Hot spots in a grounded conductor}

From a physical point of view, the solution \eqref{ck} describes the evolution of the
temperature of $M$ when its initial value distribution is known on $M$. The situation is
more difficult if $\pa M$ is not empty. We shall consider here the case of a {\it grounded} heat conductor, that is we will study the solution $h$ of  the Cauchy-Dirichlet problem
\eqref{heat1}-\eqref{heat2}.

\medskip
\noindent
{\bf Bounded conductor.}
As already seen, if $\fhi\ge 0$, \eqref{spectral}  implies  \eqref{large times}.
For an arbitrary continuous function $\fhi$, from \eqref{spectral} we can infer that, if $m$ is the first integer such that $\widehat{\fhi}(n)\not=0$ and $m+1, \dots, m+k-1$ are all the integers
such that $\la_m=\la_{m+1}=\cdots=\la_{m+k-1}$, then 
$$
\,e^{\la_m t}\, h(x,t)\to \sum_{n=m}^{m+k-1} \widehat{\fhi}(n)\,\phi_n(x) \ \mbox{ if } \  t\to\infty.
$$
Also, when $\fhi\equiv 1$, \eqref{varadhan} holds and hence 
\begin{equation}
\label{short times}
\dist(x(t), \cH_0)\to 0 \ \mbox{ as } \ t\to 0,
\end{equation}
where $\cH_0$ is the set of local (strict) maximum points of $d_\Ga$.
These informations give a rough picture of the {\it set of trajectories} of the hot spots:
$$
\cT=\bigcup_{t>0}\cH(t).
$$

\medskip

Notice in passing that, if $\Om$ is convex and  has $N$ distinct hyperplanes of symmetry, it is clear that 
$\cT$ is made of the same single point --- the intersection of the hyperplanes --- that is the hot spot {\it does not move} or is {\it stationary}. Also, it is not difficult to show (see \cite{ChS}) that the hot spot does not move if $\Om$ is invariant under an
{\it essential} group $G$ of orthogonal transformations (that is for every $x\not=0$ there is $A\in G$ such that $Ax\not=0$). Characterizing the class $\cP$ of convex domains that admit a stationary hot spot seems to be a difficult task: some partial results about convex polygons can be found in \cite{MS1,MS2} (see also \cite{MS}). There it is proved that: (i) the equilateral triangle and the parallelogram are the only polygons with $3$ or $4$ sides in $\cP$; (ii) the equilateral pentagon and the hexagons invariant under rotations of angles $\pi/3, 2\pi/3$, or $\pi$ are the only polygons with $5$ or $6$ sides {\it all} touching the inscribed circle centered at the hot spot.

\medskip

The analysis of the behavior of $\cH(t)$ for $t\to 0^+$ and $t\to\infty$ helps us to show that hot spots {\it do move} in general. 

\begin{figure}[h]
\label{fig:half-disk}
\centering
\includegraphics[scale=.8]{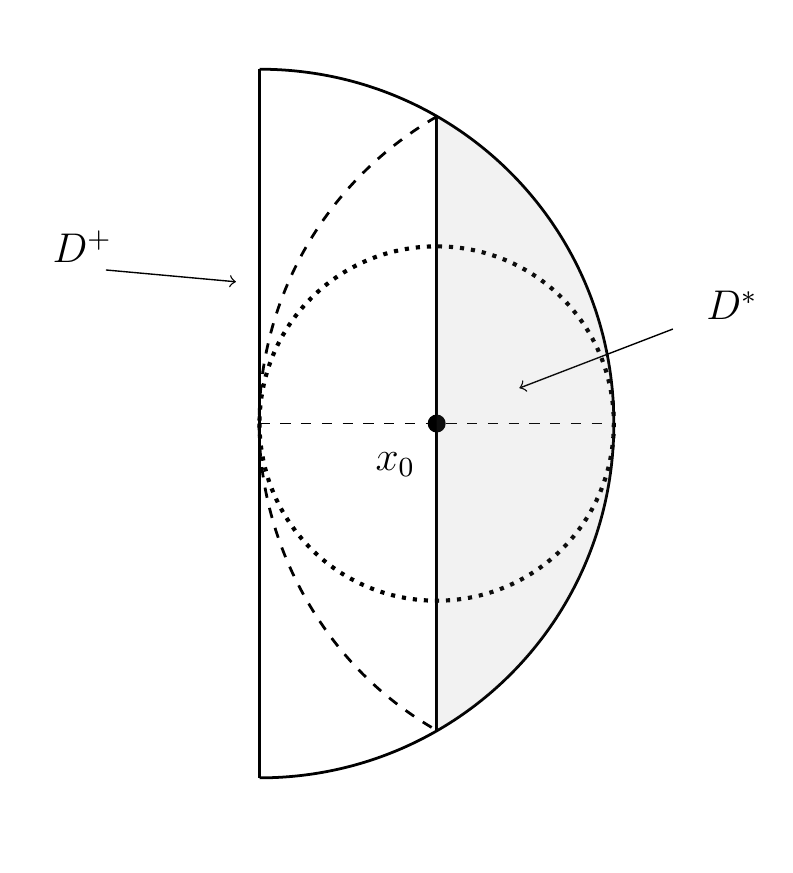}
\caption{The reflected $D^*$ is contained in $D^+$, hence $h'$ can be defined in~$D^*$.}
\end{figure}


To see this, it is enough to consider the {\it half-disk} (see Fig. 4.3)
$$
D^+=\{ (x,y)\in\RR^2: |x|<1, \ x_1>0\};
$$
being $D^+$ convex, for each $t>0$, there is a unique hot spot that, as $t\to 0^+$, tends to 
the maximum point $x_0=(1/2,0)$ of $d_\Ga$. Thus, it is enough to show that  $x_0$
is not a spatial critical point of $h(x,t)$ for some $t>0$ or, if you like, for~$\phi_1$. 
\par
This is readily
seen by {\it Alexandrov's reflection principle}. Let $D^*\!=\!\{x\in D^+: x_1>1/2\}$ and define
$$
h'(x_1,x_2,t)=h(1-x_1,x_2,t) \ \mbox{ for } \ (x_1,x_2,t)\in \ol{D^*}\times(0,\infty);
$$
$h'$ is the reflection of $h$ in the line $x_1=1/2$. We clearly have that
\begin{eqnarray*}
&(h'-h)_t=\De (h-h') \ \mbox{ in } \ D^*\times(0,\infty), \quad h'-h=0 \ \mbox{ on } \ D^*\times\{ 0\},\\
&h'\!-\!h\!>\!0 \ \mbox{ on } \ (\pa D^*\cap\pa^+ )\times(0,\infty), \quad h'\!-\!h\!=\!0 \ \mbox{ on } \
(\pa D^*\cap D^+ )\times(0,\infty).
\end{eqnarray*}
Thus, the strong maximum principle and the Hopf's boundary point lemma imply that
$$
-2\,h_{x_1}(1/2,x_2,t)=h'_{x_1}(1/2,x_2,t)-h_{x_1}(1/2,x_2,t)>0 
$$
for  $(1/2,x_2,t)\in(\pa D^*\cap D^+ )\times(0,\infty)$, and hence $x_0$ cannot be a critical point of $h$.

\medskip

The Alexandrov's principle just mentioned can also be employed to estimate the location
of a hot spot. In fact, as shown in \cite{BMS}, by the same arguments one can prove that hot
spots must belong to the subset $\core(\Om)$ of $\Om$ defined as follows. Let $\pi_\om$ be a hyperplane orthogonal to the direction $\om\in\SS^{N-1}$ and let $H^+_\om$ and $H^-_\om$ be the two half-spaces
defined by $\pi_\om$; let $\cR_\om(x)$ denote the mirror reflection of a point $x$ in $\pi_\om$.
Then, the {\it heart} \footnote{ $\core(\Om)$ has also been considered in \cite{O} under the name of {\it minimal unfolded region}.} of $\Om$ is defined by
$$
\core(\Om)=\bigcap_{\om\in\SS^{N-1}}\{H^-_\om\cap\Om : \cR_\om(H^+_\om\cap\Om)\subset\Om\}.
$$
When $\Om$ is convex, then $\core(\Om)$ is also convex and, if $\Ga$ is of class $C^1$, we are sure that its distance from $\Ga$ is positive (see \cite{Fr}). Also,  we know that $\cH(t)$ is made of only one point $x(t)$, so that 
$$
\dist(x(t),\Ga)\ge\dist(\core(\Om),\Ga).
$$
The set $\core(\Om)$ contains many notable geometric points of the set $\Om$, such as the {\it center of mass}, the {\it incenter}, the {\it circumcenter}, and others; see \cite{BM}, where further properties of
the heart of a convex body are presented. See also \cite{Sk} for related research on this issue.

\medskip

As clear from \cite{BMS}, the estimate just presented is of purely geometric nature, that is it only depends on the lack of symmetry of $\Om$ and does not depend on the particular equation we
are considering in $\Om$, as long as the equation is invariant by reflections.  
\par
A different way to estimate the location of the hot spot of a grounded convex heat conductor or the maximum point of the solution of certain elliptic equations is
based on ideas related to Alexandrov-Bakelman-Pucci's maximum principle and does take into account the information that comes from the relevant equation. For instance,  in \cite{BMS} it is proved that the maximum point $x_\infty$ of $\phi_1$ in $\ol{\Om}$ is such that
\begin{equation}
\label{bms}
\dist(x_\infty,\Ga)\ge C_N\,r_\Om\,\left(\frac{r_\Om}{\diam(\Om)}\right)^{N^2-1},
\end{equation}
where $C_N$ is a constant only depending on $N$, $r_\Om$ is the {\it inradius} of $\Om$ (the radius of a largest ball contained in $\Om$) and $\diam(\Om)$ is the {\it diameter} of $\Om$.

\medskip 

The idea of the proof of \eqref{bms} is to compare the {\it concave envelope} $f$ of $\phi_1$ --- the smallest concave function above $\phi_1$ --- and the function $g$ whose graph is the surface of the (truncated) cone based on $\Om$ and having its tip at the point $(x_\infty, \phi(x_\infty))$ (see Fig. 4.4).

\begin{figure}[h]
\label{fig:envelope}
\centering
\includegraphics[scale=.75]{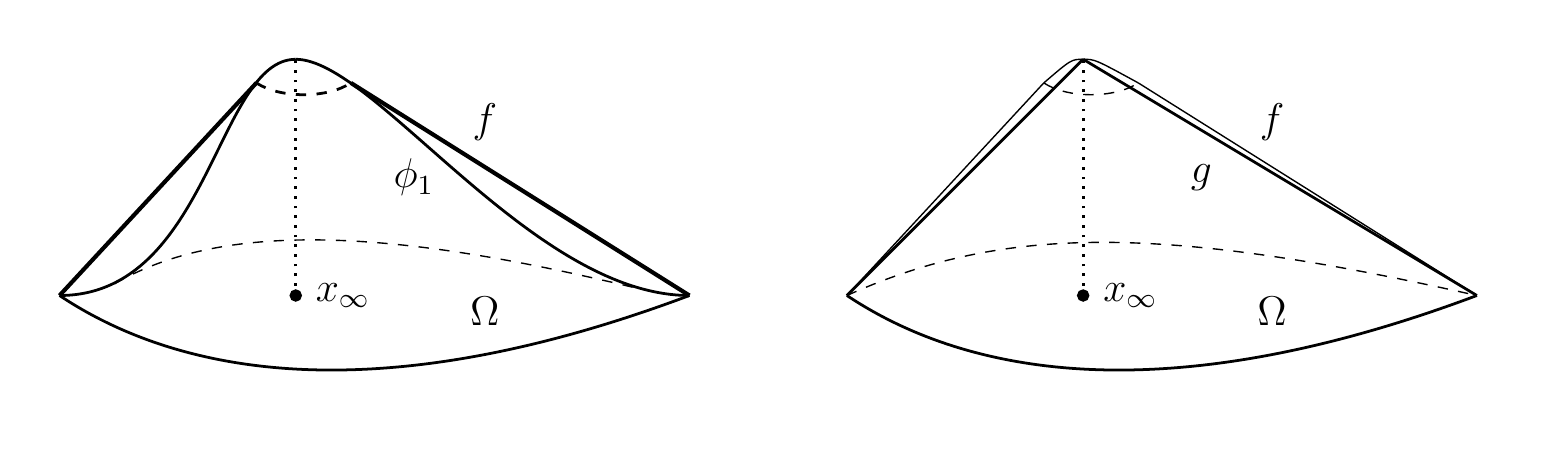}
%
\caption{The concave envelope of $\phi_1$ and the cone $g$. The dashed cap is the image $f(C)=\phi_1(C)$ of the contact set $C$.}
\end{figure}

\par
Since
$f\ge g$ and $f(x_\infty)=g(x_\infty)$, we can compare their respective {\it sub-differential images}:
\begin{eqnarray*}
&\pa f(\Om)=\bigcup\limits_{x\in\ol{\Om}}\left\{
p\in\RR^N: f(x)+p\cdot(y-x)\ge f(y) \ \mbox { for } \ y\in\ol{\Om}\right\},\\
&\pa g(\Om)=\bigcup\limits_{x\in\ol{\Om}}\left\{
p\in\RR^N: g(x)+p\cdot(y-x)\ge g(y) \ \mbox { for } \ y\in\ol{\Om}\right\};
\end{eqnarray*}
in fact, it holds that $\pa g(\Om)\subseteq\pa f(\Om)$. 


Now, $\pa g(\Om)$ has a precise geometrical meaning: it is the 
set $\phi_1(x_\infty)\,\Om^*$, that is a multiple of the {\it polar set} of $\Om$ with respect to $x_\infty$
defined by
$$
\Om^*=\{y\in\RR^N: (x-x_\infty)\cdot(y-x_\infty)\le 1 \ \mbox{ for every } \ x\in\ol{\Om}\}.
$$
The volume $|\pa f(\Om)|$ can be estimated by the formula of change of variables to obtain:
$$
\phi_1(x_\infty)^N |\Om^*|=|\pa g(\Om)|\le|\pa f(\Om)|\le\int_C |\det D^2 f|\,dx=\int_C |\det D^2 \phi_1|\,dx,
$$
where $C=\{ x\in\ol{\Om}: f(x)=\phi_1(x)\}$ is the {\it contact set}. 
Since the determinant and the trace of a matrix are the product and the sum of the eigenvalues of the matrix, by the
{\it arithmetic-geometric mean inequality}, we have that
$|\det D^2 \phi_1|\le (-\De\phi_1/N)^N$, and hence
 we can infer that
$$
|\Om^*|\le \,\int_C  \left[\frac{-\De \phi_1}{N \phi_1(x_\infty)}\right]^Ndx=
\int_C  \left[\frac{\la_1(\Om)\,\phi_1}{N \phi_1(x_\infty)}\right]^Ndx\le
\left[\frac{\la_1(\Om)}{N}\right]^N |\Om|,
$$
being $\phi_1\le\phi_1(x_\infty)$ in $\ol{\Om}$.
Finally, in order to get \eqref{bms} explicitly, one has to bound $|\Om^*|$ from below 
by the volume of the polar set of a suitable half-ball containing $\Om$, and $\la_1(\Om)$
from above by the {\it isodiametric inequality} (see \cite{BMS} for details).

The two methods we have seen so far, give estimates of how far the hot spot must be from the boundary. We now present a method, due to Grieser and Jerison  \cite{GJ}, that gives an estimate
of how far the hot spot can be from a specific point in the domain. The idea is to adapt the classical
{\it method of separation of variables} to construct a suitable
approximation $u$ of the first Dirichlet eigenfunction $\phi_1$ in a planar convex domain. Clearly, if $\Om$ were a rectangle, say $[a,b]\times[0,1]$, then that approximation would be exact: in fact
$$
u(x,y)=\phi_1(x,y)=\sin[\pi (x-a)/(b-a)].
$$
\par
If $\Om$ is not a rectangle, after some manipulations, we can suppose that 
$$
\Om=\{(x,y): a< x< b, f_1(x)<y<f_2(x)\}
$$
where, in $[a,b]$, $f_1$ is convex, $f_2$ is concave and
$$
0\le f_1\le f_2\le 1 \ \mbox{ and } \ \min_{[a,b]} f_1=0, \ \max_{[a,b]} f_2=1
$$
(see Fig. 4.5).

\begin{figure}[h]
\label{fig:long convex}
\centering
\includegraphics[scale=.8]{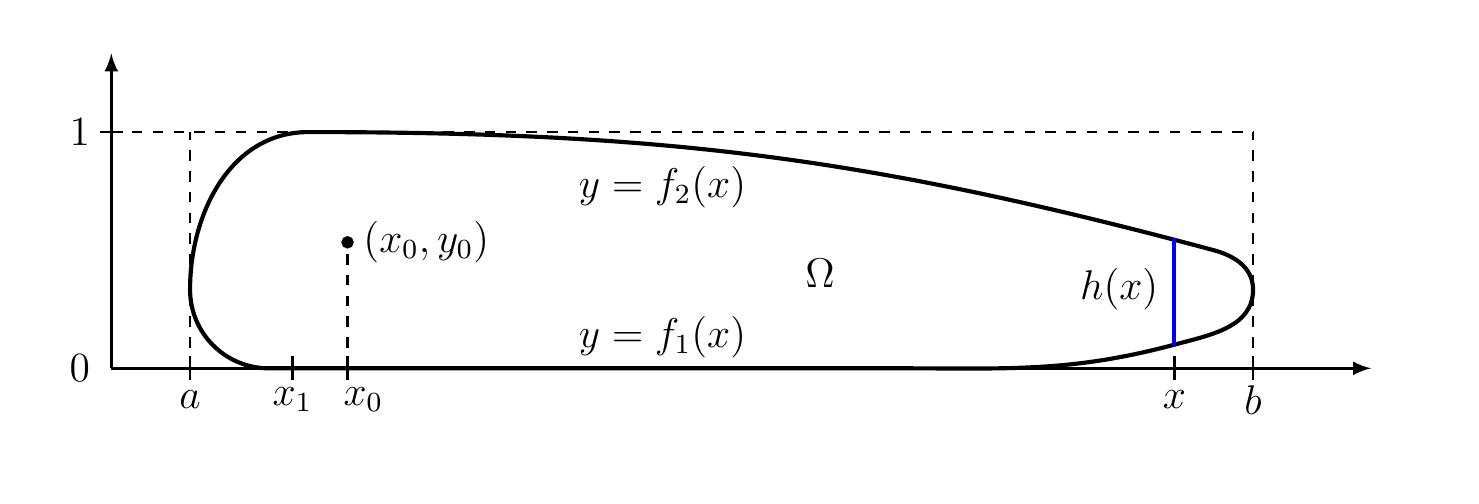}
\caption{Estimating the hot spot in the ``long'' convex set $\Om$.}
\end{figure}


The geometry of $\Om$ does not allow to find a solution by separation of variables as in the case of
the rectangle. However, one can operate ``as if'' that separation were possible.
To understand that, consider the length of the section of foot $x$, parallel to the $y$-axis, by
$$
h(x)=f_2(x)-f_1(x) \ \mbox{ for } \ a\le x\le b,
$$
and notice that, if we set 
$$
 \al(x,y)=\pi\,\frac{y-f_1(x)}{h(x)}, 
$$
the function
$$
e(x,y)=\sqrt{2/h(x)}\,\sin\al(x,y)  \ \mbox{ for } \  f_1(x)\le y\le f_2(x),
$$
satisfies for fixed $x$ the problem
$$
e_{yy}+\pi^2 e=0 \ \mbox{ in } \ (f_1(x), f_2(x)), \quad e(x,f_1(x))=e(x, f_2(x))=0
$$
--- thus, it is the first Dirichlet eigenfunction in the interval $(f_1(x), f_2(x))$, normalized in the space $L^2([f_1(x), f_2(x)])$.    The basic idea is then that $\phi_1(x,y)$ should be (and in fact it is) well approximated by its lowest Fourier mode in the $y$-direction, computed for each fixed $x$, that is by the projection of
$\phi_1$ along~$e$:
$$
\psi(x)\,e(x,y) \ \mbox { where } \ \psi(x)=\int_{f_1(x)}^{f_2(x)} \phi_1(x,\eta)\,e(x,\eta)\,d\eta.
$$ 
To simplify matters, a further approximation is needed: it turns out that $\psi$ and its first derivative can be well 
approximated by $\phi/\sqrt{2}$ and its derivative, where $\phi$ is the first  eigenfunction of the problem
$$
\phi''(x)+\left[\mu-\frac{\pi^2}{h(x)^2}\right] \phi(x)=0 \ \mbox{ for } \ a<x<b, \quad \phi(a)=\phi(b)=0.
$$
\par
Since near the maximum point $x_1$ of $\phi$, $|\phi'(x)|$ can be bounded from below by a constant
times $|x-x_1|$, the constructed chain of approximations gives that, if $(x_0,y_0)$ is the maximum point of $\phi_1$ on $\ol{\Om}$, then there is an absolute constant $C$ such that
$$
|x_1-x_0|\le C.
$$
$C$ is independent of $\Om$, but the result has clearly no content unless $b-a>C$.

\medskip
\noindent
{\bf Unbounded conductor.}
If $\Om$ is unbounded, by working with suitable barriers, one can still prove formula \eqref{varadhan} when $\fhi\equiv 1$ (see \cite{MS3,MS4}), the convergence holding uniformly on
compact subsets of $\Om$. Thus, any
hot spot $x(t)$ will again satisfy \eqref{short times}. 
\par
To the author's knowledge, \cite{JS} is the only reference in which the behavior of hot spots for large times has been studied for some grounded unbounded conductors. There, the cases of a half-space $\RR_+^N=\{x\in\RR^N: x_N>0\}$ and the exterior of a ball $B^c=\{ x\in\RR^N: |x|>1\}$ are considered.
It is shown that there is a time $T>0$ such that for $t>T$ the set $\cH(t)$ is made of only one
hot spot $x(t)=(x_1(t),\dots, x_N(t))$ and
$$
x_j(t)\to\frac{\int_{\RR^{N-1}}y_j y_N\fhi(y') dy'}{\int_{\RR^{N-1}}y_N\fhi(y')dy'}, \ 1\le j\le N-1, \quad
\frac{x_N(t)}{\sqrt{2t}}\to 1 \ \mbox{ as } \ t\to\infty,
$$
 if $\Om=\RR^N_+$, while for $\Om=B^c$, if $\fhi$ is radially symmetric, then there is a time $T>0$
 such that $\cH(t)=\{ x\in\RR^N: |x|=r(t)\}$, for $t>T$, where $r(t)$ is some smooth function of $t$
 such that 
 $$
 \limsup_{t\to\infty} r(t)=\infty.
 $$ 
Upper bounds for $\cH(t)$ are also given in \cite{JS} for the case of the exterior of a smooth bounded domain.

\subsection{Hot spots in an insulated conductor}
We conclude this survey by giving an account on the so-called {\it hot spot conjecture}
by J.~Rauch \cite{Ra}. 
This is related to the asymptotic behavior of hot spots in a {\it perfectly insulated}
heat conductor modeled by the following initial-boundary value problem:
\begin{equation}
\label{heat-neumann}
h_t=\De h \ \mbox{ in } \ \Om\times(0,\infty), \quad h=\fhi \ \mbox{ on } \ \Om\times\{ 0\},
\quad \pa_\nu u=0 \ \mbox{ on } \ \Ga\times(0,\infty).
\end{equation}
\par
Observe that, similarly to \eqref{spectral},  a spectral formula also holds for the solution of \eqref{heat-neumann}:
\begin{equation}
\label{spectral-neumann}
h(x,t)=\sum_{n=1}^\infty \widehat{\fhi}(n)\,  \psi_n(x)\, e^{-\mu_n t}, \ \mbox{ for } \ x\in\ol{\Om} \ \mbox{ and } \ t>0.
\end{equation}
Here $\{\mu_n\}_{n\in\NN}$ is the increasing sequence of Neumann eigenvalues and $\{\psi_n\}_{n\in\NN}$ is a complete orthonormal system in $L^2(\Om)$ of eigenfunctions corresponding to the $\mu_n$'s, that is $ \psi_n$ is a non-zero solution of 
\begin{equation}
\label{neumann}
\De\psi+\mu\,\psi=0 \ \mbox{ in } \ \Om, \quad \pa_\nu\psi=0 \ \mbox{ on } \ \Ga,
\end{equation}
with $\mu=\mu_n$.
The numbers $\widehat{\fhi}(n)$ are the Fourier coefficients of $\fhi$ corresponding to $\psi_n$, that is
$$
\widehat{\fhi}(n)=\int_\Om \fhi(x)\,\psi_n(x)\,dx, \ n\in\NN.
$$
Since $\mu_1=0$ and $\psi_1=1/\sqrt{|\Om|}$, we can 
infer that 
\begin{equation}
\label{spectral-neumann1}
e^{\mu_m t} \left[h(x,t)-\frac1{\sqrt{|\Om|}}\,\int_\Om\fhi\,dx\right]\to
\sum_{n=m}^{m+k-1} \widehat{\fhi}(n)\,\psi_n(x)  \ \mbox{ as } \  t\to\infty,
\end{equation}
where $m$ is the first integer such that $\widehat{\fhi}(n)\not=0$ and $m+1, \dots, m+k-1$ are all the integers
such that $\mu_m=\mu_{m+1}=\cdots=\mu_{m+k-1}$.
Thus, similarly to what happens for the case of a grounded conductor, as $t\to\infty$,
a hot spot $x(t)$ of $h$ tends to a maximum point of the function at the right-hand side of \eqref{spectral-neumann1}.
\par
Now, roughly speaking, the conjecture states that, for ``most'' initial conditions $\fhi$, the distance from 
$\Ga$ of any hot and cold spot of $h$ must tend to zero as $t\to\infty$, and hence it amounts to prove that the right-hand side of \eqref{spectral-neumann1} attains its maximum and minimum at points in $\Ga$.
\par
It should be noticed now that the quotes around the word {\it most} are justified by the fact that the conjecture does not hold for all initial conditions. In fact, as shown in \cite{BB}, if $\Om=(0,2\pi)\times(0,2\pi)\subset\RR^2$, the function defined by
$$
h(x_1,x_2, t)=-e^{-t} (\cos x_1+\cos x_2), \ (x_1,x_2)\in\Om, \ t>0,
$$
is a solution of \eqref{heat-neumann} --- with $\fhi(x_1,x_2)=-(\cos x_1+\cos x_2)$ --- that attains its maximum at $(-\pi,\pi)$ for any $t>0$. However, it turns out that in this case 
$
h(x_1,x_2, t)=-e^{-\mu_4 t} \psi_4(x_1,x_2).
$
Thus, it is wiser to rephrase the conjecture by asking whether or not the hot and cold spots tend to $\Ga$ if the coefficient $\widehat{\fhi}(2)$ of the first non-constant eigenfunction $\psi_2$ is not zero or, which is the same, whether or not maximum and minimum points of $\psi_2$ in $\ol{\Om}$ are attained only on~$\Ga$.
\par
In \cite{Ka}, a weaker version of this last statement is proved to hold for domains of the form $D\times(0,a)$, where $D\subset\RR^{N-1}$ has a boundary of class $C^{0,1}$. 
In \cite{Ka}, the conjecture has also been reformulated for convex domains. Indeed, 
we now know that it is false for fairly general domains: in \cite{BW} a planar domain with two holes is constructed, having a simple second eigenvalue and such that the corresponding eigenfunction attains its strict maximum at an interior point of the domain. It turns out that in that example the minimum point is on the boundary. Nevertheless, in \cite{BaB} it is given an example of a domain 
whose second Neumann eigenfunction attains both its maximum and minimum points at interior points. In both examples the conclusion is obtained by probabilistic methods. 
\par
Besides in \cite{Ka}, positive results on this conjecture can be found in \cite{AB,BB,BPP,Do,JN,Mi,Pa,Si}. In \cite{BB}, the conjecture is proved for planar convex domains $\Om$ with two orthogonal axis of symmetry and such that 
$$
\frac{\mbox{diam($\Om$)}}{\mbox{width($\Om$)}}>1.54.
$$
This restriction is removed in \cite{JN}. In \cite{Pa}, $\Om$ is assumed to have only one axis of symmetry, but $\psi_2$ is assumed anti-symmetric in that axis. A more general result is contained in \cite{AB}: the conjecture holds true for domains of the type
$$
\Om=\{(x_1,x_2): f_1(x_1)<x_2<f_2(x_1)\},
$$
where $f_1$ and $f_2$ have unitary Lipschitz constant. In \cite{Do}, a modified version is considered:
it holds true for general domains, if {\it vigorous maxima} are considered (see \cite{Do} for the definition). If no symmetry is assumed for a convex domain $\Om$, Y. Miyamoto \cite{Mi} has verified the conjecture when
$$
\frac{\diam(\Om)^2}{|\Om|}<1.378
$$
(for a disk, this ratio is about 1.273).

\medskip

For unbounded domains, the situation changes. For the half-space, Jimbo and Sakaguchi proved in \cite{JS} that there is a time $T$ after which the hot spot equals a point on the boundary that depends on $\fhi$. In \cite{JS}, the case of the exterior $\Om$ of a ball $\ol{B_R}$ is also considered for a radially symmetric $\fhi$. For a suitably general $\fhi$, Ishige \cite{Is1} has proved that the behavior of the hot spot is governed by the point
$$
A_\fhi=\frac{\displaystyle \int_\Om x\,\left(1+\frac{R^N}{N-1}\,|x|^{-N}\right)\,\fhi(x)\,dx}{\displaystyle
\int_\Om\fhi(x)\,dx}.
$$
If $A_\fhi\in B_R$, then $\cH(t)$ tends to the boundary point $R\,A_\fhi/|A_\fhi|$, while if
$A_\fhi\notin B_R$, then $\cH(t)$ tends to $A_\fhi$ itself. 
\par
Results concerning the behavior of hot spots for parabolic equations with a rapidly decaying potential can be found in \cite{IK1,IK2}.


\begin{thebibliography}{99}


\bibitem{Al1}
{\sc G. ~Alessandrini}, {\em Critical points of
solutions of elliptic equations in two variables}, Ann. Scuola
Norm. Super. Pisa Cl. Sci. (4) {\bf 14} (1987), 229--256. 

\bibitem{Al2} 
{\sc G.~Alessandrini}, {\em Stable determination of conductivity by boundary measurements}, Appl. Anal. {\bf 27} (1988), 153--172.

\bibitem{Al3} 
{\sc G.~Alessandrini}, {\em Singular solutions of elliptic equations and the determination of conductivity by boundary measurements}, 
J. Differential Equations {\bf 84} (1990), 252--272. 

\bibitem{AM1}
{\sc G.~Alessandrini and R.~Magnanini},
{\em The index of isolated critical points and solutions of elliptic equations in the plane}, Ann. Scuola Norm. Super. Pisa Cl. Sci. (4) {\bf 19} (1992), 567--589.

\bibitem{AM2}
{\sc G.~Alessandrini and R.~Magnanini},
{\em Elliptic equations in divergence form, geometric critical points of solutions, and Stekloff eigenfunctions}, SIAM J. Math. Anal. {\bf 25} (1994), 1259--1268.

\bibitem{AM3}
{\sc G.~Alessandrini and R.~Magnanini}, {\em Symmetry and
non--symmetry for the overdetermined Stekloff eigenvalue
problem}, Z. Angew. Math. Phys. {\bf 45} (1994), 44--52. 

\bibitem{AM4} 
{\sc G.~Alessandrini and R.~Magnanini}, {\em Symmetry and non-symmetry for the overdetermined Stekloff eigenvalue problem. II}, in Nonlinear problems in applied mathematics, 1--9, SIAM, Philadelphia, PA, 1996. 

\bibitem{AR}
{\sc G.~Alessandrini, D.~Lupo and E.~Rosset}, {\em Local behavior and geometric properties of solutions to 
degenerate quasilinear elliptic equations in the plane}, Appl. Anal. {\bf 50} (1993), 191--215.

\bibitem{AB}
{\sc R.~Atar and K.~Burdzy}, {\em On Neumann eigenfunctions in lip domains},
J. Amer. Math. Soc. {\bf 17} (2004), 243--265. 

\bibitem{BB} 
 {\sc R.~Ba\~nuelos and K.~Burdzy}, {\em On the "hot spots'' conjecture of J. Rauch}, J. Funct. Anal. {\bf 164} (1999), 1--33.
 
\bibitem{BPP} 
{\sc R.~Ba\~nuelos, M.~Pang and M.~Pascu}, {\em Brownian motion with killing and reflection and the "hot-spots'' problem}, 
Probab. Theory Related Fields {\bf 130} (2004), 56--68. 

\bibitem{BaB}
{\sc R.~Bass and K.~Burdzy}, {\em Fiber Brownian motion and the "hot spots'' problem},
Duke Math. J. {\bf 105} (2000), 25--58. 

\bibitem{BS} 
{\sc S.~Bergmann and M.~Schiffer}, Kernel Functions and Differential Equations in Mathematical Physics. Academic Press, New York, 1953.

\bibitem{Be}
{\sc L.~Bers}, {\em Function--theoretical
properties of solutions of partial differential equations of elliptic type}, Ann. Math. Stud.
{\bf 33} (1954), 69--94.

\bibitem{BN}
{\sc L.~Bers and L.~Nirenberg}, {\em On a representation theorem for linear systems with discontinuous coefficients and its applications}, in Convegno Internazionale sulle Equazioni Lineari alle Derivate Parziali, Cremonese, Roma 1955.
 
\bibitem{BiS} 
{\sc M.~Bianchini and P.~Salani}, {\em Power concavity for solutions of nonlinear elliptic problems in convex domains}, in Geometric Properties for Parabolic and Elliptic PDE’s, 35--48, Springer, Milan, 2013. 

\bibitem{BL} 
{\sc H.~J.~Brascamp and E.H.~~Lieb},  {\em On extensions of the Brunn-Minkowski and Pr\'ekopa-Leindler theorems, including inequalities for log concave functions, and with an application to the diffusion equation}, J. Funct. Anal. {\bf 22} (1976), 366--389.

\bibitem{BMS}
{\sc L.~Brasco, R.~Magnanini and P.~Salani}, {\em The location of the hot spot in a grounded convex conductor}, Indiana Univ. Math. J. {\bf 60} (2011), 633--659.

\bibitem{BM} 
{\sc L.~Brasco and R.~Magnanini}, {\em The heart of a convex body}, in Geometric properties for parabolic and elliptic PDE's, 49--66, Springer, Milan, 2013. 

\bibitem{BW} 
{\sc K.~Burdzy and W.~Werner}, {\em A counterexample to the "hot spots'' conjecture}, Ann. of Math. {\bf 149} (1999), 309--317. 

\bibitem{CS} 
{\sc L.~A.~Caffarelli and J.~Spruck}, {\em Convexity properties of solutions to some classical variational problems},
Comm. Partial Differential Equations {\bf 7} (1982), 1337--1379. 

\bibitem{CF} 
{\sc L.~A.~Caffarelli and A.~Friedman}, {\em Convexity of solutions of semilinear elliptic equations}, Duke Math. J. {\bf 52} (1985), 431--456.
 
 \bibitem{CM} 
{\sc S.~Cecchini and R.~Magnanini}, {\em Critical points of solutions of degenerate elliptic equations in the plane}, Calc. Var. Partial Differential Equations {\bf 39} (2010), 121--138. 

\bibitem{ChS} 
{\sc M.~Chamberland and D.~Siegel}, {\em Convex domains with stationary hot spots}, Math. Methods Appl. Sci. {\bf 20} (1997), 1163--1169.
 
\bibitem{CK} 
{\sc I.~Chavel and L.~Karp}, {\em Movement of hot spots in Riemannian manifolds},
J. Anal. Math. {\bf 55} (1990), 271--286. 

\bibitem{CNV} 
{\sc J.~Cheeger, A.~Naber and D.~Valtorta}, {\em Critical sets of elliptic equations}, Comm. Pure Appl. Math. {\bf 68} (2015), 173--209.

\bibitem{Do}
{\sc H.~Donnelly}, {\em Maxima of Neumann eigenfunctions},
J. Math. Phys. {\bf 49} (2008), 043506, 3 pp. 

\bibitem{EP1} 
{\sc A.~Enciso and D.~Peralta-Salas}, {\em Critical points and level sets in exterior boundary problems}, Indiana Univ. Math. J. {\bf 58} (2009), 1947--1969.

\bibitem{EP2} 
{\sc A.~Enciso and D.~Peralta-Salas}, {\em Critical points of Green's functions on complete manifolds}, J. Differential Geom. {\bf 92} (2012), 1--29.

\bibitem{EP3} 
{\sc A.~Enciso and D.~Peralta-Salas}, {\em Critical points and geometric properties of Green's functions on open surfaces}, Ann. Mat. Pura Appl. (4) {\bf 194} (2015), 881--901. 

\bibitem{EP4} 
{\sc A.~Enciso and D.~Peralta-Salas}, {\em Eigenfunctions with prescribed nodal sets}, J. Differential Geom. {\bf 101} (2015), 197--211. 

\bibitem{EJN} 
{\sc A.~Eremenko, D.~Jakobson and N.~Nadirashvili}, {\em On nodal sets and nodal domains on $\cS^2$ and $\RR^2$}, Ann. Inst. Fourier {\bf 57} (2007), 2345--2360.

\bibitem{Ev} 
{\sc L.~C.~Evans}, {\em Partial Differential Equations}, 
American Mathematical Society, Providence, RI, 1998. 

\bibitem{Fr} 
{\sc L.~E.~Fraenkel}, {\em An introduction to maximum principles and symmetry in elliptic problems}, Cambridge University Press, Cambridge, 2000.

\bibitem{Fc} 
{\sc E.~Francini}, {\em Starshapedness of level sets for solutions of nonlinear parabolic equations}, Rend. Istit. Mat. Univ. Trieste {\bf 28} (1996), 49--62.

\bibitem{Ga} 
{\sc K.~F.~Gauss}, {\em Lehrsatz}, Werke, {\bf 3}, p. 112; {\bf 8}, p.32, 1816.

\bibitem{GJ} 
{\sc D.~Grieser and D.~Jerison}, {\em The size of the first eigenfunction of a convex planar domain}, 
J. Amer. Math. Soc. {\bf 11} (1998), 41--72. 

\bibitem{Ha} 
{\sc Q.~Han}, {\em Nodal sets of harmonic functions}, Pure Appl. Math. Q. {\bf 3} (2007), 647--688. 

\bibitem{HHHN} 
{\sc R.~Hardt, M.~Hostamann-Ostenhof, T.~Hostamann-Ostenhof and N.~Nadirashvili},
{\em Critical sets of solutions to elliptic equations}, J. Differential Geom. {\bf 51} (1999), 359--373.

\bibitem{HW} 
{\sc P.~Hartman and A.~Wintner}, {\em On the local behavior of solutions of non-parabolic partial differential equations}, Amer. J. Math. {\bf 75} (1953), 449--476.
 
 \bibitem{Is1} 
{\sc K.~Ishige}, {\em Movement of hot spots on the exterior domain of a ball under the Neumann boundary condition}, J. Differential Equations {\bf 212} (2005), 394--431.
 
\bibitem{Is2} 
{\sc K.~Ishige}, {\em Movement of hot spots on the exterior domain of a ball under the Dirichlet boundary condition}, Adv. Differential Equations {\bf 12} (2007), 1135--1166.
  
\bibitem{IK1} 
{\sc K.~Ishige and Y.~Kabeya}, {\em Hot spots for the heat equation with a rapidly decaying negative potential}, Adv. Differential Equations {\bf 14} (2009), 643--662. 
 
\bibitem{IK2} 
{\sc K.~Ishige and Y.~Kabeya}, {\em Hot spots for the two dimensional heat equation with a rapidly decaying negative potential}, Discrete Contin. Dyn. Syst. Ser. S, {\bf 4} (2011), 833--849.

\bibitem{IK3} 
{\sc K.~Ishige and Y.~Kabeya}, {\em $L^p$ norms of nonnegative Schr\"odinger heat semigroup and the large time behavior of hot spots}, J. Funct. Anal. {\bf 262} (2012), 2695--2733.

\bibitem{IS} 
{\sc K.~Ishige and P.~Salani}, {\em Parabolic power concavity and parabolic boundary value problems}, Math. Ann. {\bf 358} (2014), 1091–1117.

\bibitem{Je} 
{\sc J.~L.~W.~V.~Jensen}, {\em Recherches sur la th\'eorie des \'equations},
Acta Math. {\bf 36} (1913), 181--195.

\bibitem{JaN}
{\sc D.~Jakobson and N.~Nadirashvili}, {\em Eigenfunctions with few critical points}, J. Differential Geom. {\bf 53} (1999), 177--182.

\bibitem{JNT} 
{\sc D.~Jakobson, N.~Nadirashvili and D.~Toth}, {\em Geometric properties of eigenfunctions}, Russian Math. Surveys {\bf 56} (2001), 67--88. 

\bibitem{JN}
{\sc D.~Jerison and N.~Nadirashvili}, {\em The "hot spots'' conjecture for domains with two axes of symmetry}, J. Amer. Math. Soc. {\bf 13} (2000), 741--772.

\bibitem{JS} 
{\sc S.~Jimbo and S.~Sakaguchi}, 
{\em Movement of hot spots over unbounded domains in $\RR^N$}, 
J. Math. Anal. Appl. {\bf 182} (1994), 810--835. 

\bibitem{JZ1} 
{\sc J.~Jung and S.~Zelditch}, {\em Number of nodal domains of eigenfunctions on non-positively curved surfaces with concave boundary}, Math. Ann. {\bf 364} (2016), 813--840. 

\bibitem{JZ2} 
{\sc J.~Jung and S.~Zelditch}, {\em Number of nodal domains and singular points of eigenfunctions of negatively curved surfaces with an isometric involution}, J. Differential Geom. {\bf 102} (2016), 37--66.

\bibitem{KY} 
{\sc M.~Kalka and D.~Yang}, {\em On nonpositive curvature functions on noncompact surfaces of finite topological type}, Indiana Univ. Math. J. {\bf 43} (1994), 775--804.

\bibitem{Ka} 
{\sc B.~Kawohl}, {\em Rearrangements and convexity of level sets in PDE}, Springer, Berlin, 1985.

\bibitem{Ke} 
{\sc A.~U.~Kennington}, {\em Power concavity and boundary value problems},
Indiana Univ. Math. J. {\bf 34} (1985), 687--704. 

\bibitem{Ko} 
{\sc N.~J.~Korevaar}, {\em Convex solutions to nonlinear elliptic and parabolic boundary value problems}, Indiana Univ. Math. J. {\bf 32} (1983), 603--614.

\bibitem{KL} 
{\sc N.~J.~Korevaar and J.~L.~Lewis}, 
{\em Convex solutions of certain elliptic equations have constant rank Hessians},
Arch. Ration. Mech. Anal. {\bf 97} (1987), 19--32. 

\bibitem{KP}
{\sc S.~Kantz and H.~R.~Parks}, {\em A Primer of Real Analytic Functions}, Birkh\"auser, Basel, 2002.

\bibitem{LV} 
{\sc O. Lehto and K. Virtanen}, {\em Quasiconformal Mappings in the Plane}, Springer, Berlin, 1973.

\bibitem{LT} 
{\sc P.~Li and L.~F.~Tam}, {\em Symmetric Green's functions on complete manifolds}, Amer. J. Math. {\bf 109} (1987), 1129--1154. 

\bibitem{Lu} 
{\sc F.~Lucas}, {\em Propri\'et\'es g\'eom\'etriques des fractions rationelles}, Paris Comptes Rendus {\bf 78} (1874), 271--274.

\bibitem{MS}
{\sc R.~Magnanini and S.~Sakaguchi},  {\em The spatial critical points not moving along the heat flow},
J. Anal. Math. {\bf 71} (1997), 237--261.

\bibitem{MS1} 
{\sc R.~Magnanini and S.~Sakaguchi}, {\em On heat conductors with a stationary hot spot}, Ann. Mat. Pura Appl. (4) 183 (2004), 1--23.
 
 \bibitem{MS2} 
{\sc R.~Magnanini and S.~Sakaguchi}, {\em Polygonal heat conductors with a stationary hot spot}, J. Anal. Math. {\bf 105} (2008), 1--18.

\bibitem{MS3}
{\sc R.~Magnanini and S.~Sakaguchi}, {\em Interaction between nonlinear diffusion and geometry of domain}, J. Differential Equations {\bf 252} (2012), 236--257.

\bibitem{MS4} 
{\sc R.~Magnanini and S.~Sakaguchi}, {\em Matzoh ball soup revisited: the boundary regularity issue},  Math. Methods Appl. Sci. {\bf 36} (2013), 2023--2032.

\bibitem{Ma} 
{\sc M.~Marden}, {\em The Geometry of the Zeros of a Polynomial in a Complex Variable}, American Mathematical Society, New York, N. Y., 1949.

\bibitem{Mi}
{\sc Y.~Miyamoto}, {\em The "hot spots'' conjecture for a certain class of planar convex domains}, J. Math. Phys. {\bf 50} (2009), 103530, 7 pp. 

\bibitem{Mo}
{\sc M.~Morse}, {\em Relations between the
critical points of a real function of n independent variables}, Trans. Amer. Math. Soc. {\bf 27} (1925), 345--396.

\bibitem{MC} 
{\sc M.~Morse and S.~S.~Cairns}, 
{\em Critical point theory in global analysis and differential topology: An introduction},
Academic Press, New York-London 1969.

\bibitem{O}
{\sc J.~O'Hara}, {\em Minimal unfolded regions of a convex hull and parallel bodies},
preprint (2012) arXiv:1205.0662v2.

\bibitem{Pa}
{\sc M.~Pascu}, {\em Scaling coupling of reflecting Brownian motions and the hot spots problem}, 
Trans. Amer. Math. Soc. {\bf 354} (2002), 4681--4702. 

\bibitem{Pe} 
{\sc D.~Peralta-Salas}, {\em private communication}, (2016).

\bibitem{Pu} 
{\sc C. Pucci}, {\em An angle's maximum principle for the gradient of solutions of
elliptic equations},  Boll. Unione Mat. Ital. {\bf 1} (1987), 135--139.

\bibitem{Ra}
{\sc J.~Rauch}, {\em Five problems: an introduction to the qualitative theory of partial differential equations}, in Partial differential equations and related topics, 355--369, Springer, Berlin, 1975. 

\bibitem{Ro} 
{\sc E.~H.~Rothe}, {\em A relation between the type numbers of a critical point
and the index of the corresponding field of gradient vectors}, Math. Nachr. {\bf 4} (1950-51), 12--27.

\bibitem{Sa} 
{\sc S. Sakaguchi}, {\em Critical points of solutions to
the obstacle problem in the plane},  Ann. Sc. Norm. Super. Pisa Cl. Sci. (4) {\bf 21} (1994), 157--173.

\bibitem{Sa2} 
{\sc S. Sakaguchi}, {\em private communication}, (2008).

\bibitem{Sk}
{\sc S.~Sakata}, {\em Movement of centers with respect to various potentials}, Trans. Amer. Math. Soc. {\bf 367} (2015), 8347--8381.
 
\bibitem{Sl1}
{\sc P.~Salani}, {\em Starshapedness of level sets of solutions to elliptic PDEs},  Appl. Anal. {\bf 84} (2005), 1185--1197.

\bibitem{Sl2} 
{\sc P.~Salani}, {\em Combination and mean width rearrangements of solutions of elliptic equations in convex sets}, Ann. Inst. H. Poincar\'e Analyse Non Lin\'eaire {\bf 32} (2015), 763--783.

\bibitem{Si}
{\sc B.~Siudeja}, {\em Hot spots conjecture for a class of acute triangles},
Math. Z. {\bf 280} (2015), 783--806. 

\bibitem{To}
{\sc J-C.~Tougeron}, {\em Id\'eaux de fonctions differentiables}, Springer, Berlin-New York, 1972. 

\bibitem{Va} 
{\sc S.~R.~S.~Varadhan}, {\em On the behavior
of the fundamental solution of the heat equation with
variable coefficients}, Comm. Pure Appl. Math. {\bf 20} (1967), 431--455.

\bibitem{Ve}
{\sc I. N. Vekua},
{\em Generalized Analytic Functions}, Pergamon Press, Oxford, 1962.

\bibitem{Wa} 
{\sc J. L. Walsh}, {\em The Location of
Critical Points of Analytic and Harmonic Functions}, 
American Mathematical Society, New York, NY, 1950.

\bibitem{Wh} 
{\sc H. Whitney}, {\em A function not constant on a connected set of critical points}, Duke Math. J. {\bf 1} (1935), 514--517. 

\bibitem{Ya}
{\sc S.~T.~Yau}, {\em Problem section, Seminar on Differential Geometry}, 
Ann. of Math. Stud. {\bf 102} (1982) 669--706. 

\bibitem{Ya2}
{\sc S.~T.~Yau}, 
{\em Selected expository works of Shing-Tung Yau with commentary. Vol. I-II}, 
International Press, Somerville, MA; Higher Education Press, Beijing, 2014. 

\bibitem{Ze} 
{\sc S.~Zelditch}, {\em private communication}, (2016).

\end{thebibliography}
\end{document}